\documentstyle[myart,amsfonts,12pt]{article}
\normalbaselines
\font\tenbm=cmmib10
\font\sevenbm=cmmib7
\font\fivebm=cmmib5
\newfam\bmfam
\textfont\bmfam=\tenbm \scriptfont\bmfam=\sevenbm
\scriptscriptfont\bmfam=\fivebm
{\count0=\number\bmfam \multiply\count0 by "100
\def\defbgreek#1#2#3{{\count1=\count0 \advance\count1 by "#2#3
  \global\mathchardef#1=\count1 }}
\defbgreek\balpha  0B \defbgreek\brho       1A
\defbgreek\bbeta   0C \defbgreek\bsigma     1B
\defbgreek\bgamma  0D \defbgreek\btau       1C
\defbgreek\bdelta  0E \defbgreek\bupsilon   1D
\defbgreek\bepsilon0F \defbgreek\bphi       1E
\defbgreek\bzeta   10 \defbgreek\bchi       1F
\defbgreek\bmeta   11 \defbgreek\bpsi       20
\defbgreek\btheta  12 \defbgreek\bomega     21
\defbgreek\biota   13 \defbgreek\bvarepsilon22
\defbgreek\bkappa  14 \defbgreek\bvartheta  23
\defbgreek\blambda 15 \defbgreek\bvarpi     24
\defbgreek\bmu     16 \defbgreek\bvarrho    25
\defbgreek\bnu     17 \defbgreek\bvarsigma  26
        \defbgreek\bxi     18 \defbgreek\bvarphi    27
\defbgreek\bpi     19}
\newtheorem{opred}{Definition}[section]
\input tcilatex.tex

\begin{document}

\title{Geometry without Topology as a New Conception of Geometry}
\author{Yuri A. Rylov}
\date{Institute for Problems in Mechanics, Russian Academy of Sciences, \\
101-1, Vernadskii Ave., Moscow, 117526, Russia. \\
email: rylov@ipmnet.ru\\
Web site: {$http://rsfq1.physics.sunysb.edu/\symbol{126}rylov/yrylov.htm$}.\\
or mirror Web site: {$http://194.190.131.172/\symbol{126}rylov/yrylov.htm$}.}
\maketitle

\begin{abstract}
A geometric conception is a method of a geometry construction. The
Riemannian geometric conception and a new T-geometric one are considered.
T-geometry is built only on the basis of information included in the metric
(distance between two points). Such geometric concepts as dimension,
manifold, metric tensor, curve are fundamental in the Riemannian conception
of geometry, and they are derivative in the T-geometric one. T-geometry is
the simplest geometric conception (essentially only finite point sets are
investigated) and simultaneously it is the most general one. It is
insensitive to the space continuity and has a new property -- nondegeneracy.
Fitting the T-geometry metric with the metric tensor of Riemannian geometry,
one can compare geometries, constructed on the basis of different
conceptions. The comparison shows that along with similarity (the same
system of geodesics, the same metric) there is a difference. There is an
absolute parallelism in T-geometry, but it is absent in the Riemannian
geometry. In T-geometry any space region is isometrically embeddable in the
space, whereas in Riemannian geometry only convex region is isometrically
embeddable. T-geometric conception appears to be more consistent logically,
than the Riemannian one.
\end{abstract}

\newpage

\section{Introduction}

Conception of geometry (geometric conception) is a method (a set of
principles), which is used for construction of geometry. The proper
Euclidean\footnote{%
We use the term ''Euclidean geometry'' as a collective concept with respect
to terms ''proper Euclidean geometry'' and ''pseudoeuclidean geometry''. In
the first case the eigenvalues of the metric tensor matrix have similar
signs, in the second case they have different signs. The same interrelation
take place between terms "Riemannian geometry", "proper Riemannian geometry"
and "pseudo-Riemannian geometry".} geometry can be constructed on the basis
of different geometric conceptions. For instance, one can use the Euclidean
axiomatic conception (Euclidean axioms), or the Riemannian conception of
geometry (dimension, manifold, metric tensor, curve). One can use metric
conception of geometry (topological space, metric, curve). In any case one
obtains the same proper Euclidean geometry. From point of view of this
geometry it is of no importance which of possible geometric conceptions is
used for the geometry construction. It means that the category of geometry
conception is metageometric.

However, if we are going to generalize (to modify) the Euclidean geometry,
it appears to be very important which of many possible geometric conceptions
is used for the generalization. The point is that the generalization is some
modification of original (fundamental) statements of geometry in the scope
of the same geometric conception. As far as fundamental statements are
different in different geometric conceptions, one is forced to modify
different statements, that leads naturally to different results.

If one uses the Euclidean geometric conception, which contains only axioms
and no numerical characteristics, the only possible modification consists in
changing some axioms by other ones. In this case some new geometries appear
which hardly may be considered to be a generalization of the Euclidean
geometry. They are rather its different modifications.

Some fundamental statements of the Riemannian geometric conception contain
numerical characteristics, as far as one sets the dimension $n$ and metric
tensor $g_{ik}, \quad i,k=1,2,\ldots n$, consisting of several functions of
one point, i.e. of one argument $x=\{x^i\},\;\; i=1,2,\ldots n$. Varying $n$
and $g_{ik}$, one obtains a class of Riemannian geometries, where each
geometry is labelled by several functions of one point.

Recently a new geometric conception of the Euclidean geometry construction
was suggested \cite{R90,R92}. The Euclidean geometry appears to be
formulated in terms only of metric $\rho $, setting distance between all
pairs of points of the space. Such a geometric conception is the most
general in the sense, that all information on geometry is concentrated in
one function of two points. It is evident that one function of two points
contains more information, than several functions of one point (it is
supposed that the set of points is continual). At some choice of the point
set $\Omega $, where the metric and geometry are set, the $n$-dimensional
Euclidean geometry appears. At another choice of the metric another
generalized geometry appears on the same set $\Omega $. This geometry will
be referred to as tubular geometry, or briefly T-geometry. All things being
equal, the set of all T-geometries appears to be more powerful, than the set
of all Riemannian geometries. This conception will be referred to as
T-geometric conception, although the term ''metric conception of geometry''
fits more.

The point is that this term has been occupied. By the metric (or generalized
Riemannian) geometry \cite{T59,ABN86,BGP92} is meant usually a geometry,
constructed on the basis of the metric geometric conception, whose
fundamental statements are topology and metric, i.e. the metric is set not
on an arbitrary set of points, but on the topological space, where, in
particular, concepts of continuity and of a curve are defined.

What actually is happen is that the metric geometric conception contains
excess of fundamental statements. This excess appears as follows. Let us
imagine that some conception $A$ of Euclidean geometry contains some set of
independent fundamental statements $a$. Let $b$ be some set of corollaries
of the fundamental statements $a$. Let us consider now the set $a\cup b$ as
a set of fundamental statements of a geometric conception. It is another
conception $A^{\prime }$ of Euclidean geometry. Its fundamental statements $%
a\cup b$ are not independent. Now one can obtain the conception $A$, or some
other geometric conception, depending on how the fundamental statements $%
a\cup b$ are used. Now obtaining generalized geometries, one may not vary
the fundamental statements independently. To avoid contradiction, one is to
take into account mutual dependence of fundamental statements.

If we know nothing on mutual dependence of fundamental statements $a\cup b$,
the geometric conception may appear to be eclectic. We risk to obtain
contradictions, or artificial constraints on the generalized geometries
obtained. In the case of the metric conception of the proper Euclidean
geometry the statements on properties of the topological space and those of
the curve are corollaries of metrical statements. They may be removed
completely from the set of fundamental statements of the conception.

However, there are problems, connected with the fact that we have some
preconceptions on what is the geometry, in general. In particular, it is a
common practice to consider that the concept of the curve is an attribute of
any geometry, which is used for description of the real space (or
space-time). It is incorrect, and manifests itself, in particular, in
imposition of some unjustified constraints (triangle inequality) on metric,
which make the difficult situation. These preconceptions have a metalogic
character. They are connected with association properties of human thinking.
Overcoming of these preconceptions needs a serious analysis.

A cause for writing this paper is a situation, arising after appearance and
discussion of papers on T-geometry \cite{R90,R92}, which mean essentially a
construction of a new geometric conception. Such a situation took place in
the second half of XIXth century, when the non-Euclidean geometries
appeared, and the most part of mathematical community considered sceptically
applications of the Riemannian (and non-Euclidean) geometry to the real
space geometry. Appearance of Riemannian geometries meant appearance of a
new geometric conception. The reason of sceptical relation of the
mathematical community to Riemannian geometry has not been analyzed up till
now, although it was described in literature \cite{K37}.

Appearance of a more general conception of geometry is important for
applications of geometry. Geometry is a ground of the space-time model, and
appearance of a new more general geometric conception poses the question as
to whether the microcosm space-time geometry has been chosen optimally. If
the existing space-time geometry is not optimal, it must be revised. The
space-time geometry revision is to be accompanied by a revision of basic
statements of physics as a science founded on the space-time model. For
instance, appearance of Riemannian geometries and realization of the fact,
that a new conception of geometry appears together with their appearance,
has lead finally to a revision of the space-time conception and to creation
of the general relativity theory.

Until appearance of T-geometries there was only one uniform isotropic
geometry suitable for the space-time description. This is Minkowski
geometry. An alternative to the Minkowski geometry to be anywhere reasonable
did not exist. After realization of the fact that non-degenerate geometries
(T-geometries) are as good as degenerate (Riemannian) geometries, a class of
geometries suitable for description of uniform isotropic space-time appears.
This class includes the Minkowski geometry. The uniform isotropic geometries
of this class are labelled by a function of one argument. Geometries of the
class differ in a value and character of nondegeneracy. All geometries of
this class except for Minkowski geometry appear to be nondegenerate. In the
nondegenerate geometry a motion of free particles appears to initially
stochastic, whereas in the degenerate geometry it initially deterministic.
It is well known, that motion of microparticles (electrons, protons, etc.)
is stochastic. It seems incorrect to choose such a space-time model, where
the microparticle motion is deterministic, and thereafter to introduce
additional hypotheses (principles of quantum mechanics), providing
stochasticity of microparticle motion. It is more reasonable to choose at
once such a space-time geometry which provides the microparticle motion
stochasticity. It is desirable to choose from the class of uniform
nondegenerate geometries precisely that geometry, which agrees optimally
with experimental data. If the complete agreement with experiment appears to
be impossible, one can add supplementary hypotheses, as it is made in
quantum mechanics. In any case the space-time geometry is to be chosen
optimally. The choice of the Minkowski geometry as a space-time model for
microcosm is not optimal certainly. A use of the Minkowski geometry as a
space-time model for microcosm is explained by absence of alternative (i.e.
essentially by a use of the Riemannian conception of geometry).

Thus, after appearance of a new conception of geometry and appearance of an
alternative to the Minkowski geometry a revision of the space-time model is
a logical necessity. This revision must be carried out independently of that
whether the new version of the space-time model explains all quantum
effects, or only part of them. In the last case one should add some
hypotheses, explaining that part of experimental data, which are not
explained by the revised space-time model. In any case one should use the
most suitable space-time geometry among all possible ones.

Let us note that this conclusion does not agree with viewpoint of most of
physicists, dealing with relativistic quantum theory. Many of them suppose
that any revision of the existing space-time model is justified only in the
case, if it explains at least one of experiments which cannot be explained
by the existing theory. We agree with such a position, provided the existing
theory modification does not concern principles of quantum theory and
space-time model. At appearance of a more general conception of geometry one
is forced to choose an optimal geometry independently of whether the new
model solves all problems, or only a part of them. Another viewpoint, when
one suggests either to solve all problems by means of a revision of the
space-time geometry, or, if it appears to be impossible, to abandon from
revision at all and to use certainly nonoptimal geometry, seems to be too
extremistic.

Now results of application of nondegenerate geometry for the space-time
description seem to be rather optimistic, because one succeeded to choose
such a nondegenerate geometry, containing the quantum constant $\hbar$ as a
parameter, that statistical description of stochastic particle motion in
this space-time coincides with the quantum description in terms of
Schr\"odinger equation in the conventional space-time model \cite{R91,R1995}%
. Further development of the conception will show whether explanation of
relativistic quantum effects is possible.

In the present paper a new geometric conception, based on the concept of
distance and {\sl only distance} is considered. In general, the idea of the
geometry construction on the basis of the distance is natural and not new.
The geometric conception, where the distance (metric) is a basic concept,
 is natural to be referred to as metric
conception of geometry. Usually the term ''metric geometry'' is used for a
geometry, constructed on the base of the metric space.

\begin{opred}
\label{d0} The metric space $M=\{\rho ,\Omega \}$ is the set $\Omega $ of
points $P\in \Omega $, equipped by the metric $\rho $, setting on $\Omega
\times \Omega $
\begin{equation}
\rho :\quad \Omega \times \Omega \rightarrow D_{+}\subset {\Bbb R},\qquad
D_{+}=[0,\infty ),  \label{a1.1}
\end{equation}
\begin{equation}
\rho (P,P)=0,\qquad \rho (P,Q)=\rho (Q,P),\qquad \forall P,Q\in \Omega
\label{a1.2}
\end{equation}
\begin{equation}
\rho (P,Q)=0,\quad \mbox{if and only if}{\rm \;\;\;\;}P=Q,\qquad \forall
P,Q\in \Omega   \label{a1.3}
\end{equation}
\begin{equation}
\rho (P,Q)+\rho (Q,R)\geq \rho (P,R),\qquad \forall P,Q,R\in \Omega
\label{a1.4}
\end{equation}
\end{opred}

There is a generalization of metric geometry, known as distance geometry
\cite{B53}, which differs from the metric geometry in absence of constraint (%
\ref{a1.4}). The main problem of metric geometric conception is a construction
of geometric objects, i.e. different sets of points in the metric space. For
instance, to construct such a geometric object as the shortest, one is
forced to introduce the concept of a curve as a continuous mapping of a
segment of real axis on the space.
\begin{equation}
L:\;I\rightarrow \Omega ,\qquad I=[0,1]\subset {\Bbb R},  \label{b5.5}
\end{equation}
The shortest, passing through points $P$ and $Q$, is defined as a curve
segment of the shortest length. On one hand,
introduction of the concept of a curve means a rejection from the pure
metric conception of geometry, as far as one is forced to introduce
concepts, which do not defined via metric. On the other hand, if the concept
of a curve is not introduced, it is not clear how to build such geometric
objects which are analogs of Euclidean straight and plane. Without
introduction of these objects the metric geometry looks as a very pure
(slightly informative) geometry. Such a geometry cannot be used as a model
of the real space-time.

Essentially the problem of constructing a pure metric conception of
geometry is set as follows. Is it possible to construct on the basis of only
metric such a geometry which were as informative as the Euclidean geometry?
In other words, is it possible to construct the Euclidean geometry, setting
in some way the metric on $\Omega \times \Omega $, where $\Omega $ is a
properly chosen set of points? More concretely this problem is formulated as
follows.

Let $\rho _{{\rm E}}$ be the metric of $n$-dimensional proper Euclidean
space on $\Omega\times\Omega$. Is it possible on the base of information,
contained in
$\rho _{{\rm E}}$ to reconstruct the Euclidean geometry, i.e. to determine
the dimension $n$, to introduce rectilinear coordinate system and metric
tensor in it, to construct $k$-dimensional planes $k=1,2,\ldots n$ and to
test whether the reconstructed geometry is proper Euclidean? If yes, and
information, contained in metric is sufficient for construction of
proper Euclidean
geometry, the used prescriptions can be used for construction of a geometry
with other metric. As a result each metric $\rho$ corresponds to some metric
geometry $T_\rho$, constructed {\sl on the base of the metric and only metric%
}. Any such a geometry $T_\rho$ is not less as informative as the proper Euclidean
one in the sense, that any geometric object in proper Euclidean geometry
corresponds to a geometric object in the metric geometry $T_\rho$,
constructed according to the same prescriptions, as it is built in the
proper Euclidean geometry. This geometric object may appear to bear little
resemblance to its Euclidean analog. Besides, due to symmetry of the
Euclidean space (presence of a motion group) different geometric objects in $%
T_\rho$ may have the same Euclidean analog. For instance, in the Euclidean
geometry any two different points $P$ and $Q$, lying on the Euclidean
straight $L$, determine this straight. In metric geometry $T_\rho$ analogs
of the Euclidean straight ${\cal T}_{PQ}$, ${\cal T}_{P_1Q_1}$, $P_1,Q_1\in
{\cal T}_{PQ}$, determined by different pairs $P,Q$ and $P_1,Q_1$, are
different, in general, if the metric does not satisfy the condition (\ref
{a1.4}).

There exists a positive solution of the stated problem, i.e. amount of
information, contained in the metric, is sufficient for constructing the
metric geometry which is not less informative, than the Euclidean one.
Corresponding theorem has been proved \cite{R01}.

Apparently K.~Menger \cite{M28} succeeded to approach most closely to the
positive solution of the mentioned problem, but he failed to solve it
completely. The reason of his failure is some delusion, which may be
qualified as ''associative prejudice''. An overcoming of this prejudice
results a construction of new geometric conception, where all information on
geometry is contained in metric. The new conception generates a class of
T-geometries, which may be considered to be a generalization of conventional
metric geometry on the base of metric space $M=\{\rho ,\Omega \}$. Formally
this generalization is approached at the expense of reduction of number of
fundamental concepts, i.e. concepts necessary for the geometry construction
and at the expense of elimination of constraints (\ref{a1.3}), (\ref{a1.4}),
imposed on metric. Besides instead of metric $\rho $ one uses the quantity $%
\sigma =\frac{1}{2}\rho ^{2}$, known as world function \cite{S60}. The world
function is supposed to be real. It means that the metric $\rho $ may be
either nonnegative, or pure imaginary quantity. This extends capacities of
geometry. Now one can consider the Minkowski geometry as a special case of
T-geometry and use the T-geometry as a space-time geometry. The concept of
the curve (\ref{b5.5}) is not used at the construction of geometry, i.e. it
is not a fundamental concept, although as the geometry construction has
been completed completed nothing prevents from introduction of the curve
by means of the mapping (\ref{b5.5}).

But the curve $L$ appears not to be an attribute of geometry. It is some
additional object external with respect to geometry. A corollary of this is
an appearance of a new geometry property, which is referred to as
nondegeneracy. Euclidean and Riemannian geometries have no nondegeneracy.
They are degenerate geometries.
Associative prejudice is an delusion, appearing, when properties of one
object are attributed by mistake to another object. Let us illustrate this
in a simple example, which is perceived now as a grotesque. It is known that
ancient Egyptians believed that all rivers flow towards the North. This
delusion seems now to be nonsense. But many years ago it had weighty
foundation. The ancient Egyptians lived on a vast flat plane and knew only
one river the Nile, which flew exactly towards the North and had no
tributaries on the Egyptian territory. The North direction was a preferred
direction for ancient Egyptians who observed motion of heavenly bodies
regularly. It was direction toward the fixed North star. They did not
connect direction of the river flow with the plane slope, as we do now. They
connected the direction of the river flow with the preferred spatial
direction towards the North. We are interested now what kind of mistake made
ancient Egyptians, believing that all rivers flow towards the North, and how
could they to overcome their delusion.

Their delusion was not a logical mistake, because the logic has no relation
to this mistake. The delusion was connected with associative property of
human thinking, when the property $A$ is attributed to the object $B$ on the
basis that in all known cases the property $A$ accompanies the object $B$.
Such an association may be correct or not. If it is erroneous, as in the
given case, it is very difficult to discover the mistake. At any rate it is
difficult to discover the mistake by means of logic, because such
associations appear before the logical analysis, and the subsequent logical
analysis is carried out on the basis of the existing associations. Let us
imagine that in the course of a voyage an ancient Egyptian scientist arrived
the Tigris, which is the nearest to Egypt river. He discovers a water stream
which flows, first, not outright and, second, not towards the North. Does he
discover his delusion? Most likely not. At any rate not at once. He starts
to think that the water stream, flowing before him, is not a river. A ground
for such a conclusion is his initial belief that ''real'' river is to flow,
first, directly and, second, towards the North. Besides, the Nile was very
important in the life of ancient Egyptians, and they were often apt to
idolize the Nile. The delusion about direction of the river flow can be
overcame only after that, when one has discovered sufficiently many
different rivers, flowing towards different directions, and the proper
analysis of this circumstance has been carried out.

Thus, to overcome the associating delusion, it is not sufficient to present
another object $B$, which has not the property $A$, because one may doubt of
whether the presented object is to be classified really as the object $B$.
Another attendant circumstances are also possible.

If the established association between the object and its property is
erroneous, one can say on associative delusion or on associative prejudice.
The usual method of overcoming the associative prejudices is a consideration
of wider set of phenomena, where the established association between the
property $A$ and the object $B$ may appear to be violated, and the
associative prejudice is discovered.

The associative prejudices are very stable. It is very difficult to overcome
them, when they have been established, because they cannot be disproved
logically. On the other hand, fixing incorrect correlations between objects
of real world, the associative prejudices point out a wrong way for
investigations.

Associative prejudices are known in history of science. For instance, the
known statement of the Ptolemaic doctrine that the Earth is placed in the
centre of universe, and heaven bodies rotate around it, is an example of the
associative prejudice. In this case the property of being a centre of a
planetary system is attributed to the Earth, whereas such a centre is the
Sun. Overcoming of this prejudice was long and difficult, because in
contrast to prejudice of ancient Egyptians it can be disproved neither
logically, nor experimentally.

Another example of associative prejudice is the popular in XIXth century
opinion that the Cartesian coordinate system is an attribute of geometry.
This view point appeared, when the analytic geometry was discovered, and the
Cartesian coordinate system became to serve as a tool at description of
geometric objects of Euclidean geometry. Using analytic description of
Euclidean geometry, many mathematicians of XIXth century applied Cartesian
coordinates almost always and were inclined to believe that the Cartesian
coordinates are an attribute of any geometry at all. On the other hand,
non-uniform (Riemannian) geometry cannot be constructed in the Cartesian
coordinate system. Any attempt of writing the Riemannian geometry metric
tensor in a Cartesian coordinates turns non-uniform (Riemannian) geometry to
uniform (i.e. Euclidean) geometry. In other words, the Cartesian coordinate
system discriminates any non-uniform geometry. It is known \cite{K37} that
mathematicians of XIXth century were biased against consideration of the
Riemannian geometry as a really existing geometry. It seems that this
scepsis in the relation of Riemannian geometry is connected with the
associative prejudice, when the Cartesian coordinate system is considered to
be an attribute of any geometry. As the coordinate system appears to be a
way of the geometry description, but not its attribute, the scepsis
disappears.

Now the viewpoint that the concept of the curve (\ref{b5.5}) is a
fundamental concept (i.e. it is applied at construction of any geometry)
holds much favor. This viewpoint is based on the circumstance that the curve
is used at construction of all known (Riemannian and metric) geometries.
Such a viewpoint is an associative prejudice (of the kind as the statement
of ancient Egyptians that all rivers flow towards the North). To prove this,
it is sufficient to construct a sufficiently informative geometry without
using the concept of the curve (\ref{b5.5}). Such a geometry (T-geometry) has
been constructed \cite{R90}. Constructing the new conception of geometry,
its author did not think that he did not use the concept of the curve and
overcame some prejudice. The point is that the metric $\rho (x,y)$,
considered to be a function of two variable points $x$ and $y$, contains
much more information, than the metric tensor $g_{ik}(x)$, which is several
functions of one variable point $x$. The author believed that information
contained in metric is sufficient for constructing geometry, and he wants to
construct a geometry on the base of only this information. It is possible,
provided the concept of the curve is ignored. He did not suspect that he
overcame the associative preconception on fundamental role of the curve and,
hence, created a new conception of geometry. All this became clear well
later at realization and discussion of the obtained results.

In the second section the T-geometric technique is described, and one shows
that the Euclidean geometry can be formulated in terms of only metric. The
method of the geometric objects, constructed in T-geometry, is described in the
third section. The fourth section is devoted to the convexity problem. In
the fifth and sixth sections one compares solutions of the parallelism
problem in Riemannian and tubular geometries.

\section{$\protect\sigma $-space and T-geometry}

T-geometry is constructed on $\sigma $-space $V=\{\sigma ,\Omega \}$, which
is obtained from the metric space after removal of constraints (\ref{a1.3}),
(\ref{a1.4}) and introduction of the world function $\sigma $
\begin{equation}
\sigma (P,Q)\equiv \frac 12\rho ^2(P,Q),\qquad \forall P,Q\in \Omega .
\label{a1.8}
\end{equation}
instead of the metric $\rho$:

\begin{opred}
\label{d3.1.1} $\sigma $-space $V=\{\sigma ,\Omega \}$ is nonempty set $%
\Omega $ of points $P$ with given on $\Omega \times \Omega $ real function $%
\sigma $
\begin{equation}
\sigma :\quad \Omega \times \Omega \rightarrow {\Bbb R},\qquad \sigma
(P,P)=0,\qquad \sigma (P,Q)=\sigma (Q,P)\qquad \forall P,Q\in \Omega .
\label{a2.1}
\end{equation}
\end{opred}

The function $\sigma $ is known as the world function \cite{S60}, or $\sigma
$-function. The metric $\rho $ may be introduced in $\sigma $-space by means
of the relation (\ref{a1.8}). If $\sigma $ is positive, the metric $\rho $
is also positive, but if $\sigma $ is negative, the metric is imaginary.

\begin{opred}
\label{d3.1.1аe}. Nonempty point set $\Omega ^{\prime }\subset \Omega $ of $%
\sigma $-space $V=\{\sigma ,\Omega \}$ with the world function $\sigma
^{\prime }=\sigma |_{\Omega ^{\prime }\times \Omega ^{\prime }}$, which is a
contraction $\sigma $ on $\Omega ^{\prime }\times \Omega ^{\prime }$, is
called $\sigma $-subspace $V^{\prime }=\{\sigma ^{\prime },\Omega ^{\prime
}\}$ of $\sigma $-space $V=\{\sigma ,\Omega \}$.
\end{opred}

Further the world function $\sigma ^\prime = \sigma |_{\Omega ^{\prime
}\times\Omega ^\prime }$, which is a contraction of $\sigma $ will be
denoted as $\sigma $. Any $\sigma$-subspace of $\sigma$-space is a $\sigma$%
-space.

\begin{opred}
\label{d3.1.1ba}. $\sigma $-space $V=\{\sigma ,\Omega \}$ is called
isometrically embeddable in $\sigma $-space $V^{\prime }=\{\sigma ^{\prime
},\Omega ^{\prime }\}$, if there exists such a monomorphism $f:\Omega
\rightarrow \Omega ^{\prime }$, that $\sigma (P,Q)=\sigma ^{\prime
}(f(P),f(Q))$,\quad $\forall P,\forall Q\in \Omega ,\quad f(P),f(Q)\in
\Omega ^{\prime }$,
\end{opred}

Any $\sigma $-subspace $V^{\prime }$ of $\sigma $-space $V=\{\sigma ,\Omega
\}$ is isometrically embeddable in it.

\begin{opred}
\label{d3.1.1b}. Two $\sigma $-spaces $V=\{\sigma ,\Omega \}$ and $V^{\prime
}=\{\sigma ^{\prime },\Omega ^{\prime }\}$ are called to be isometric
(equivalent), if $V$ is isometrically embeddable in $V^{\prime }$, and $%
V^{\prime }$ is isometrically embeddable in $V$.
\end{opred}

\begin{opred}
\label{d2.2b} The $\sigma $-space $M=\{\rho ,\Omega \}$ is called a finite $%
\sigma $-space, if the set $\Omega $ contains a finite number of points.
\end{opred}

\begin{opred}
\label{d3.1.1bc}. The $\sigma $-subspace $M_{n}({\cal P}^{n})=\{\sigma ,%
{\cal P}^{n}\}$of the $\sigma $-space $V=\{\sigma ,\Omega \}$, consisting of
$n+1$ points ${\cal P}^{n}=\left\{ P_{0},P_{1},...,P_{n}\right\} $ is called
the $n$th order $\sigma $-subspace .
\end{opred}

The T-geometry is a set of all propositions on properties of $\sigma $%
-subspaces of $\sigma $-space $V=\{\sigma ,\Omega \}$. Presentation of
T-geometry is produced on the language, containing only references to $%
\sigma $-function and constituents of $\sigma $-space, i.e. to its $\sigma $%
-subspaces.

\begin{opred}
\label{d3} A description is called $\sigma $-immanent, if it does not
contain any references to objects or concepts other, than finite subspaces
of the metric space and its world function (metric).
\end{opred}

$\sigma $-immanence of description provides independence of the description
on the method of description. In this sense the $\sigma $-immanence of a
description in T-geometry reminds the concept of covariance in Riemannian
geometry. Covariance of some relation in Riemannian geometry means that the
considered relation is valid in all coordinate systems and, hence, describes
only the properties of the Riemannian geometry in itself. Covariant
description provides cutting-off from the coordinate system properties,
considering the relation in all coordinate systems at once. The $\sigma $%
-immanence provides truncation from the methods of description by absence of
a reference to objects, which do not relate to geometry itself (coordinate
system, concept of curve, dimension).

The basic elements of T-geometry are finite $\sigma $-subspaces $M_n({\cal P}%
^n)$, i.e. finite sets
\begin{equation}
{\cal P}^n=\{P_0,P_1,\ldots ,P_n\}\subset \Omega  \label{a1.9}
\end{equation}

The main characteristic of the finite $\sigma $-subspace $M_{n}({\cal P}^{n})
$ is its length $|M\left( {\cal P}^{n}\right) |$

\begin{opred}
\label{d1.13} The squared length $\left| M\left( {\cal P}^{n}\right) \right|
^{2}$ of the $n$th order $\sigma $-subspace $M\left( {\cal P}^{n}\right)
\subset \Omega $ of the $\sigma $-space $V=\left\{ \sigma ,\Omega \right\} $
is the real number.
\[
\left| M\left( {\cal P}^{n}\right) \right| ^{2}=\left( n!S_{n}({\cal P}%
^{n})\right) ^{2}=F_{n}\left( {\cal P}^{n}\right)
\]
where $S_{n}({\cal P}^{n})$ is the volume of the $(n+1)$-edr, whose
vertices are placed
at points ${\cal P}^{n}\equiv \{P_{0},P_{1},\ldots P_{n}\}\subset \Omega $,
defined by means of relations
\begin{equation}
F_n:\quad \Omega ^{n+1}\rightarrow {\Bbb R},\qquad \Omega
^{n+1}=\bigotimes\limits_{k=1}^{n+1}\Omega ,\qquad n=1,2,\ldots  \label{a1.5}
\end{equation}
\begin{equation}
F_n\left( {\cal P}^n\right) =\det ||\left( {\bf P}_0{\bf P}_i.{\bf P}_0{\bf P%
}_k\right) ||,\qquad P_0,P_i,P_k\in \Omega ,\qquad i,k=1,2,...n  \label{a1.6}
\end{equation}
\begin{eqnarray}
\left( {\bf P}_0{\bf P}_i.{\bf P}_0{\bf P}_k\right) &\equiv &\Gamma \left(
P_0,P_i,P_k\right) \equiv \sigma \left( P_0,P_i\right) +\sigma \left(
P_0,P_k\right) -\sigma \left( P_i,P_k\right) ,  \label{a1.7} \\
i,k &=&1,2,...n,  \nonumber
\end{eqnarray}
\end{opred}
where the function $\sigma $ is defined via metric $\rho $ by the relation (%
\ref{a1.8}) and ${\cal P}^n$ denotes $n+1$ points (\ref{a1.9}).

The meaning of the written relations is as follows. In the special case,
when the $\sigma $-space is Euclidean space and its $\sigma $-function
coincides with $\sigma $-function of Euclidean space, any two points $P_0,P_1
$ determine the vector ${\bf P}_0{\bf P}_1$, and the relation (\ref{a1.7})
is a $\sigma $-immanent expression for the scalar product $\left( {\bf P}_0%
{\bf P}_i.{\bf P}_0{\bf P}_k\right)$ of two vectors. Then the relation (\ref
{a1.6}) is the Gram's determinant for $n$ vectors ${\bf P}_0{\bf P}_i,\quad
i=1,2,\ldots n$, and $S_n({\cal P}^n)$ is the Euclidean volume of the $(n+1)$%
-edr with vertices at the points ${\cal P}^n$.

The idea of constructing the T-geometry is very simple. All relations of
proper Euclidean geometry are written in the $\sigma $-immanent form and
declared
to be valid for any $\sigma $-function. This results that any relation of
proper Euclidean geometry corresponds some relation of T-geometry. It is
important
that in the relations, declared to be relations of T-geometry, only the
properties (\ref{a1.8}) were used. The special properties of the Euclidean $%
\sigma $-function are not to be taken into account. The metric part of these
relations was formulated and proved by K.~Menger \cite{M28}. Let us present
this result in our designations in the form of the theorem
\begin{theorem}
\label{t1} The $\sigma $-space $V=\{\sigma ,\Omega \}$ is isometrically
embeddable in $n$-dimensional proper Euclidean space $E_n$, if and only
if any $(n+2)$th order $\sigma $-subspace $M({\cal P}^{n+2})\subset \Omega $
is isometrically embeddable in $E_n$.
\end{theorem}

Unfortunately, the formulation of this theorem is not $\sigma $-immanent, as
far as it contains a reference to $n$-dimensional Euclidean space $E_{n}$
which is not defined $\sigma $-immanently. A more constructive version of
the $\sigma $-space Euclideaness conditions is formulated in the form

\noindent I.
\begin{equation}
\exists {\cal P}^{n}\subset \Omega ,\qquad F_{n}({\cal P}^{n})\neq 0,\qquad
F_{n+1}(\Omega ^{n+2})=0,  \label{a3.4}
\end{equation}
II.
\begin{equation}
\sigma (P,Q)={\frac{1}{2}}\sum_{i,k=1}^{n}g^{ik}({\cal P}^{n})[x_{i}\left(
P\right) -x_{i}\left( Q\right) ][x_{k}\left( P\right) -x_{k}\left( Q\right)
],\qquad \forall P,Q\in \Omega ,  \label{a3.5}
\end{equation}
where the quantities $x_{i}\left( P\right) $, $x_{i}\left( Q\right) $ are
defined by the relations
\begin{equation}
x_{i}\left( P\right) =\left( {\bf P}_{0}{\bf P}_{i}.{\bf P}_{0}{\bf P}%
\right) ,\qquad x_{i}\left( Q\right) =\left( {\bf P}_{0}{\bf P}_{i}.{\bf P}%
_{0}{\bf Q}\right) ,\qquad i=1,2,...n  \label{a3.5a}
\end{equation}
The contravariant components $g^{ik}({\cal P}^{n}),$ $(i,k=1,2,\ldots n)$ of
metric tensor are defined by its covariant components $g_{ik}({\cal P}^{n}),$
$(i,k=1,2,\ldots n)$ by means of relations
\begin{equation}
\sum_{k=1}^{n}g_{ik}({\cal P}^{n})g^{kl}({\cal P}^{n})=\delta
_{i}^{l},\qquad i,l=1,2,\ldots n  \label{a3.11}
\end{equation}
where
\begin{equation}
g_{ik}({\cal P}^{n})=\Gamma (P_{0},P_{i},P_{k}),\qquad i,k=1,2,\ldots n
\label{a3.9}
\end{equation}
III.\quad The relations
\begin{equation}
\Gamma (P_{0},P_{i},P)=x_{i},\qquad x_{i}\in {\Bbb R},\qquad i=1,2,\ldots n,
\label{a3.12}
\end{equation}
considered to be equations for determination of $P\in \Omega $, have always
one and only one solution.

\noindent IIIa. The relations (\ref{a3.12}), considered to be equations
for determination of $P\in \Omega $, have
always not more than one solution.

\begin{remark}
\label{r3} The condition (\ref{a3.4}) is a corollary of the condition
(\ref{a3.5}). It is formulated in the form of a special condition, in order
that a determination of dimension were separated from determination of
coordinate system.
\end{remark}

The condition I determines the space dimension. The condition II describes $%
\sigma $-immanently the scalar product properties of the proper Euclidean space.
 Setting $n+1$ points $%
{\cal P}^{n}$, satisfying the condition I, one determines $n$-dimensional
basis of vectors in Euclidean space. Relations (\ref{a3.9}), (\ref{a3.11})
determine covariant and contravariant components of the metric tensor, and
the relations (\ref{a3.5a}) determine covariant coordinates of points $P$
and $Q$ at this basis. The relation (\ref{a3.5}) determines the expression
for $\sigma $-function for two arbitrary points in terms of coordinates of
these points. Finally, the condition III describes continuity of the set $%
\Omega $ and a possibility of the manifold construction on it. Necessity of
conditions I -- III for Euclideaness of $\sigma $-space is evident. One can
prove their sufficiency \cite{R01}. The connection of conditions I -- III
with the Euclideaness of the $\sigma $-space can be formulated in the form
of a theorem.
\begin{theorem}
\label{c2}The $\sigma $-space $V=\{\sigma ,\Omega \}$ is the $n$-dimensional
Euclidean  space, if and only if $\sigma $-immanent
conditions I -- III are fulfilled.
\end{theorem}
\begin{remark}
\label{r2} For the $\sigma$-space were proper Euclidean, the eigenvalues
of the matrix $g_{ik}({\cal P}^n),\quad i,k=1,2,\ldots n$ must have the same
sign, otherwise it is pseudoeuclidean.
\end{remark}
The theorem states that it is sufficient to know metric (world function) to
construct Euclidean geometry. The information, contained in concepts of
topological space and curve, which are used in metric geometry, appears to
be excess.

Proof of this theorem can be found in \cite{R01}. A similar theorem for
another (but close) necessary and sufficient conditions has been proved in
ref. \cite{R90}. Here we show only constructive character of conditions I --
III for proper Euclidean space. It means that starting from an abstract $%
\sigma $-space, satisfying conditions I -- III, one can determine dimension $%
n$ and construct a rectilinear coordinate system with conventional description
of the proper Euclidean space in it. One construct sequentially straight,
two-dimensional plane, etc...up to $n$-dimensional plane coincide with the
set $\Omega $. To construct all these objects, one needs to develop
technique of T-geometry.

\begin{opred}
\label{d.2.6} The finite $\sigma $-space $M_{n}({\cal P}^{n})=\{\sigma ,%
{\cal P}^{n}\}$ is called oriented $\overrightarrow{M_{n}({\cal P}^{n})}$,
if the order of its points ${\cal P}^{n}=\{P_{0},P_{1},\ldots P_{n}\}$ is
fixed.
\end{opred}

\begin{opred}
\label{d3.1.2b}. The $n$th order multivector $m_{n}$ is the mapping
\begin{equation}
m_{n}:\qquad I_{n}\rightarrow \Omega ,\qquad I_{n}\equiv \left\{
0,1,...,n\right\}   \label{a2.1dc}
\end{equation}
\end{opred}
The set $I_{n}$ has a natural ordering, which generates an ordering of
images $m_{n}(k)\in \Omega $ of points $k\in I_{n}$.
The ordered list of images of points in $I_{n}$ has one-to-one connection
with the multivector and may be used as the multivector identificator.
Different versions of the point list will be used for writing the $n$th
order multivector identificator:
\[
\overrightarrow{P_{0}P_{1}...P_{n}}\equiv {\bf P}_{0}{\bf P}_{1}...{\bf P}%
_{n}\equiv \overrightarrow{{\cal P}^{n}}
\]
Originals of points $P_{k}$ in $I_{n}$ are determined by the order of the
point $P_{k}$ in the list of identificator. Index of the point $P_{k}$ has
nothing to do with the original of $P_{k}$. Further we shall use
identificator. $\overrightarrow{P_{0}P_{1}...P_{n}}$ of the multivector
instead of the multivector. In this sense the $n$th order multivector $%
\overrightarrow{P_{0}P_{1}...P_{n}}$ in the $\sigma $-space $V=\{\sigma
,\Omega \}$ may be defined as the ordered set $\{P_{l}\},\quad l=0,1,\ldots n
$ of $n+1$ points $P_{0},P_{1},...,P_{n}$, belonging to the $\sigma $-space $%
V$. The point $P_{0}$ is the origin of the multivector $\overrightarrow{%
P_{0}P_{1}...P_{n}}$. Image $m_{n}\left( I_{n}\right) $ of the set $I_{n}$
contains $k$ points ($k\leq n+1).$ The set of all $n$th order multivectors $%
m_{n}$ constitutes the set $\Omega
^{n+1}=\bigotimes\limits_{k=1}^{n+1}\Omega $, and any multivector $%
\overrightarrow{{\cal P}^{n}}\in \Omega ^{n+1}$.

\begin{opred}
\label{d3.1.6b}. The scalar $\sigma $-product $(\overrightarrow{{\cal P}^{n}}%
.\overrightarrow{{\cal Q}^{n}})$ of two $n$th order multivectors $%
\overrightarrow{{\cal P}^{n}}$ and $\overrightarrow{{\cal Q}^{n}}$ is the
real number
\begin{equation}
(\overrightarrow{{\cal P}^{n}}.\overrightarrow{{\cal Q}^{n}})=\det \Vert (%
{\bf P}_{0}{\bf P}_{i}.{\bf Q}_{0}{\bf Q}_{k})\Vert
,\,\,\,\,\,\,\,\,\,\,\,i,k=1,2,...n,\qquad \overrightarrow{{\cal P}^{n}},%
\overrightarrow{{\cal Q}^{n}}\in \Omega ^{n+1}  \label{a2.8}
\end{equation}
\begin{eqnarray}
({\bf P}_{0}{\bf P}_{i}.{\bf Q}_{0}{\bf Q}_{k})\equiv \sigma
(P_{0},Q_{i})+\sigma (Q_{0},P_{k})-\sigma (P_{0},Q_{0})-\sigma (P_{i},Q_{k}),
\label{a2.8.1} \\
P_{0},P_{i},Q_{0},Q_{k}\in \Omega   \nonumber
\end{eqnarray}
\end{opred}

\begin{opred}
\label{d3.1.6c}. The length $|\overrightarrow{{\cal P}^{n}}|$ of the $n$th
order multivector $\overrightarrow{{\cal P}^{n}}$ is the number
\begin{equation}
|\overrightarrow{{\cal P}^{n}}|=\left\{
\begin{array}{c}
\mid \sqrt{(\overrightarrow{{\cal P}^{n}}.\overrightarrow{{\cal P}^{n}})}%
\mid ,\quad (\overrightarrow{{\cal P}^{n}}.\overrightarrow{{\cal P}^{n}}%
)\geq 0 \\
i\mid \sqrt{(\overrightarrow{{\cal P}^{n}}.\overrightarrow{{\cal P}^{n}})}%
\mid ,\quad (\overrightarrow{{\cal P}^{n}}.\overrightarrow{{\cal P}^{n}})<0
\end{array}
\right. \qquad \overrightarrow{{\cal P}^{n}}\in \Omega ^{n+1}  \label{c1.10}
\end{equation}
\end{opred}

In the case, when multivector does not contain similar points, it coincides
with the oriented finite $\sigma $-subspace, and it is a constituent of $%
\sigma $-space. In the case, when at least two points of multivector
coincide, the multivector length vanishes, and the multivector is considered
to be null multivector. The null multivector is not a finite $\sigma $%
-subspace, but a use of null multivectors assists in creation of a more
simple technique. In the case of manipulation with numbers, written in
Arabic numerals (where zero is present) is simpler, than the same
manipulation with numbers, written in Roman numerals (where zero is absent).
Something like that takes place in the case of multivectors. Essentially, the
multivectors are basic objects of T-geometry. As to continual geometric
objects, which are analogs of planes, sphere ellipsoid, etc., they are
constructed by means of skeleton-envelope method (see next section) with
multivectors, or finite $\sigma $-subspaces used as skeletons. As a
consequence the T-geometry is presented $\sigma $-immanently, i.e. without
reference to objects, external with respect to $\sigma $-space.

\begin{opred}
\label{d3.1.5c}. Two $n$th order multivectors $\overrightarrow{{\cal P}^{n}}$
$\overrightarrow{{\cal Q}^{n}}$ are collinear $\overrightarrow{{\cal P}^{n}}%
\parallel \overrightarrow{{\cal Q}^{n}}$, if
\begin{equation}
(\overrightarrow{{\cal P}^{n}}.\overrightarrow{{\cal Q}^{n}})^{2}=|%
\overrightarrow{{\cal P}^{n}}|^{2}\cdot |\overrightarrow{{\cal Q}^{n}}|^{2}
\label{a2.10}
\end{equation}
\end{opred}

\begin{opred}
\label{d3.1.5e}. Two collinear $n$th order multivectors $\overrightarrow{%
{\cal P}^{n}}$ and $\overrightarrow{{\cal Q}^{n}}$ are similarly oriented $%
\overrightarrow{{\cal P}^{n}}\uparrow \uparrow \overrightarrow{{\cal Q}^{n}}$
(parallel), if
\begin{equation}
(\overrightarrow{{\cal P}^{n}}.\overrightarrow{{\cal Q}^{n}})=|%
\overrightarrow{{\cal P}^{n}}|\cdot |\overrightarrow{{\cal Q}^{n}}|
\label{a2.11}
\end{equation}
They have opposite orientation $\overrightarrow{{\cal P}^{n}}\uparrow
\downarrow \overrightarrow{{\cal Q}^{n}}$ (antiparallel), if
\begin{equation}
(\overrightarrow{{\cal P}^{n}}.\overrightarrow{{\cal Q}^{n}})=-|%
\overrightarrow{{\cal P}^{n}}|\cdot |\overrightarrow{{\cal Q}^{n}}|
\label{a2.12}
\end{equation}
\end{opred}

Vector ${\bf P}_0{\bf P}_1=\overrightarrow{{\cal P}^1}$ is the first order
multivector.

\begin{opred}
\label{d1.14} $n$th order $\sigma $-subspace $M\left( {\cal P}^{n}\right) $
of nonzero length $\;\left| M\left( {\cal P}^{n}\right) \right|
^{2}=F_{n}\left( {\cal P}^{n}\right) \neq 0$ determines the set of points $%
{\cal T}$ $\left( {\cal P}^{n}\right) $, called $n$th order tube by means of
relation
\begin{equation}
{\cal T}\left( {\cal P}^{n}\right) \equiv {\cal T}_{{\cal P}^{n}}=\left\{
P_{n+1}|F_{n+1}\left( {\cal P}^{n+1}\right) =0\right\} ,\qquad P_{i}\in
\Omega ,\qquad i=0,1\ldots n+1,  \label{b1.3}
\end{equation}
where the function $F_{n}$ is defined by the relations (\ref{a1.5}) -- (\ref
{a1.7})
\end{opred}

In arbitrary T-geometry the $n$th order tube is an analog of $n$-dimensional
properly Euclidean plane.

\begin{opred}
\label{d3.1.9}. Section ${\cal S}_{n;P}$ of the tube ${\cal T}({\cal P}^{n})$
at the point $P\in {\cal T}({\cal P}^{n})$ is the set ${\cal S}_{n;P}({\cal T%
}({\cal P}^{n}))$ of points, belonging to the tube ${\cal T}({\cal P}^{n})$
\begin{equation}
{\cal S}_{n;P}({\cal T}({\cal P}^{n}))=\{P^{\prime }\mid
\bigwedge_{l=0}^{l=n}\sigma (P_{l},P^{\prime })=\sigma (P_{l},P)\},\qquad
P\in {\cal T}({\cal P}^{n})\qquad P^{\prime }\in \Omega .  \label{a2.38}
\end{equation}
\end{opred}

Let us note that ${\cal S}_{n;P}({\cal T}({\cal P}^n))\subset {\cal T}({\cal %
P}^n)$, because $P\in {\cal T}({\cal P}^n)$. Indeed, whether the point $P$
belongs to ${\cal T}({\cal P}^n)$ depends only on values of $n+1$ quantities
$\sigma (P_l,P),\;\;l=0,1,...n$. In accordance with (\ref{a2.38}) these
quantities are the same for both points $P$ and $P^{\prime }$. Hence, any
running point $P^{\prime }\in {\cal T}({\cal P}^n)$, if $P\in {\cal T}({\cal %
P}^n)$.

In the proper Euclidean space the $n$th order tube is $n$-dimensional plane,
containing points ${\cal P}^{n}$, and its section ${\cal S}_{n;P}({\cal T}(%
{\cal P}^{n}))$ at the point $P$ consists of one point $P$.

Now we can construct the proper Euclidean space and rectilinear coordinate
system in it on the basis of only $\sigma $-function. Let it is known
that the $\sigma $-space $V=\{\sigma ,\Omega \}$ is the proper Euclidean
space, but its dimension is not known. To determine the dimension $n$, let us
take two different points $P_0,P_1\in \Omega ,\;\;F_1({\cal P}^1)=2\sigma
(P_0,P_1)\neq 0$.

1. Let us construct the first order tube ${\cal T}\left( {\cal P}^1\right) $%
. If ${\cal T}\left( {\cal P}^1\right) =\Omega $, then dimension of the $%
\sigma $-space $V$ \ $n=1$. If $\Omega \backslash {\cal T}\left( {\cal P}%
^1\right) \neq \emptyset ,\;\;\exists P_2\in \Omega ,\;\;P_2\notin {\cal T}%
\left( {\cal P}^1\right) ,$ and hence, $F_2({\cal P}^2)\neq 0$.

2. Let us construct the second order tube ${\cal T}\left( {\cal P}^2\right) $%
. If ${\cal T}\left( {\cal P}^2\right) =\Omega $, then \ $n=2$, otherwise$%
\;\;\exists P_3\in \Omega ,\;\;P_3\notin {\cal T}\left( {\cal P}^2\right) ,$
and hence, $F_3({\cal P}^3)\neq 0$.

3. Let us construct the third order tube ${\cal T}\left( {\cal P}^3\right) $%
. If ${\cal T}\left( {\cal P}^3\right) =\Omega $, then $n=3$, otherwise$%
\;\;\exists P_4\in \Omega ,\;\;P_4\notin {\cal T}\left( {\cal P}^3\right) ,$
and hence, $F_4({\cal P}^4)\neq 0$.

4. Etc.

Continuing this process, one determines such $n+1$ points ${\cal P}^n$, that
the condition ${\cal T}\left( {\cal P}^n\right) =\Omega $ and, hence,
conditions (\ref{a3.4}) are fulfilled.

Then by means of relations
\begin{equation}
x_i\left( P\right) =\Gamma (P_0,P_i,P),\qquad i=1,2,\ldots n,  \label{a3.13}
\end{equation}
one attributes covariant coordinates $x\left( P\right) =\left\{
x_i(P)\right\} ,\;\;i=1,2,\ldots n$ to $\forall P\in \Omega $. Let $%
x=x\left( P\right) \in {\Bbb R}^n$ and $x^{\prime }=x\left( P^{\prime
}\right) \in {\Bbb R}^n.$ Substituting $\Gamma (P_0,P_i,P)=x$ and $\Gamma
(P_0,P_i,P^{\prime })=x_i^{\prime }$ in (\ref{a3.5}), one obtains the
conventional expression for the world function of the Euclidean space in the
rectilinear coordinate system
\begin{equation}
\sigma (P,P^{\prime })=\sigma _E(x,x^{\prime })={\frac 12}%
\sum_{i,k=1}^ng^{ik}({\cal P}^n)\left( x_i-x_i^{\prime }\right) \left(
x_k-x_k^{\prime }\right)  \label{a3.14}
\end{equation}
where $g^{ik}({\cal P}^n)$, defined by relations (\ref{a3.9}) and (\ref
{a3.11}), is the contravariant metric tensor in this coordinate system.

Condition III of the theorem states that the mapping
\[
x:\;\Omega \rightarrow {\Bbb R}^n
\]
described by the relation (\ref{a3.13}) is a bijection, i.e. $\forall y\in
{\Bbb R}^n$ there exists such one and only one point $Q\in \Omega ,$ that $%
y=x\left( Q\right) $.

Thus, on the base of the world function, given on abstract set $\Omega
\times \Omega $, one can determine the dimension $n$ of the Euclidean space,
construct rectilinear coordinate system with the metric tensor $g_{ik}({\cal %
P}^n)=\Gamma (P_0,P_i,P_k),\qquad i,k=1,2,\ldots n$ and describe all
geometrical objects which are determined in terms of coordinates. The
Euclidean space and Euclidean geometry is described in terms and only in
terms of world function (metric).

Conditions I -- III, formulated in the $\sigma $-immanent form admit one to
construct the proper Euclidean space, using only information, contained in
world function. $\sigma $-immanence of the formulation admits one to state
that information, contained in the world function, is sufficient for
construction of any T-geometry. Substitution of condition III by the
condition IIIa leads to a reduction of constraints. At the fulfillment of
conditions I,II,IIIa the $\sigma $-space appears to be isometrically
embeddable in $n$-dimensional Euclidean space. It may be piecewise
continuous, or even discrete. Such a $\sigma $-space can be obtained,
removing arbitrary number of points from $n$-dimensional Euclidean space.

\section{Skeleton-envelope method of geometric objects construction}

\begin{opred}
\label{d1.7} Geometric object ${\cal O}$ is some $\sigma $-subspace of $%
\sigma $-space.
\end{opred}

In T-geometry a geometric object ${\cal O}$ is described by means of
skeleton-envelope method. It means that any geometric object ${\cal O}$ is
considered to be a set of intersections and joins of  elementary geometric
objects (EGO).

\begin{opred}
Elementary geometric object ${\cal E}\subset \Omega $ is a set of zeros of
the envelope function
\begin{equation}
f_{{\cal P}^{n}}:\qquad \Omega \rightarrow {\Bbb R},\qquad {\cal P}%
^{n}\equiv \left\{ P_{0},P_{1},...P_{n}\right\} \subset \Omega   \label{b1.4}
\end{equation}
i.e..
\begin{equation}
{\cal E}={\cal E}_{f}\left( {\cal P}^{n}\right) =\left\{ R|f_{{\cal P}%
^{n}}\left( R\right) =0\right\}   \label{b1.5}
\end{equation}
\label{dd3.1}
\end{opred}

The finite set ${\cal P}^{n}\subset \Omega $ of parameters of the envelope
function $f_{{\cal P}^n}$ is skeleton of elementary geometric object (EGO)
${\cal E}\subset \Omega $. The set ${\cal E}\subset \Omega $ of points
forming EGO is called the envelope of its skeleton ${\cal P}^{n}$. For
continuous T-geometry the envelope ${\cal E}$ is usually a continual set of
points. The envelope function $f_{{\cal P}^{n}}$, determining EGO is a
function of the running point $R\in \Omega $ and of parameters ${\cal P}%
^{n}\subset \Omega $. The envelope function $f_{{\cal P}^{n}}$ is supposed
to be an algebraic function of $s$ arguments $w=\left\{
w_{1},w_{2},...w_{s}\right\} $, $s=(n+2)(n+1)/2$. Each of arguments $%
w_{k}=\sigma \left( Q_{k},L_{k}\right) $ is a $\sigma $-function of two
arguments $Q_{k},L_{k}\in \left\{ R,{\cal P}^{n}\right\} $, either belonging
to skeleton ${\cal P}^{n}$, or coinciding with the running point $R$.

Let us consider examples of some simplest EGOs.
\begin{equation}
{\cal S}(P_{0},P_{1})=\left\{ R|f_{P_{0}P_{1}}\left( R\right) =0\right\}
,\qquad f_{P_{0}P_{1}}\left( R\right) =\sqrt{2\sigma \left(
P_{0},P_{1}\right) }-\sqrt{2\sigma \left( P_{0},R\right) }  \label{b1.6}
\end{equation}
is a sphere, passing through the point $P_{1}$ and having its center at the
point $P_{0}$. Ellipsoid ${\cal EL}$, passing through the point $P_{2}$ and
having the focuses at points $P_{0},P_{1}$ $\left( P_{0}\neq P_{1}\right) $
is described by the relation
\begin{equation}
{\cal EL}(P_{0},P_{1},P_{2})=\left\{ R|f_{P_{0}P_{1}P_{2}}\left( R\right)
=0\right\} ,  \label{b1.7}
\end{equation}
where the envelope function $f_{P_{0}P_{1}P_{2}}\left( R\right) $ is defined
by the equation.
\begin{equation}
f_{P_{0}P_{1}P_{2}}\left( R\right) =\sqrt{2\sigma \left( P_{0},P_{2}\right) }%
+\sqrt{2\sigma \left( P_{1},P_{2}\right) }-\sqrt{2\sigma \left(
P_{0},R\right) }-\sqrt{2\sigma \left( P_{1},R\right) }  \label{b1.8}
\end{equation}
If focuses $P_{0},P_{1}$ coincide $\left( P_{0}=P_{1}\right) $, the
ellipsoid ${\cal EL}(P_{0},P_{1},P_{2})$ degenerates into a sphere ${\cal S}%
(P_{0},P_{2})$. If the points $P_{1},P_{2}$ coincide $\left(
P_{1}=P_{2}\right) $, the ellipsoid ${\cal EL}(P_{0},P_{1},P_{2})$
degenerates into a segment of a straight line ${\cal T}_{[P_{0}P_{1}]}$
between the points $P_{0},P_{1}$.
\begin{equation}
{\cal T}_{[P_{0}P_{1}]}={\cal EL}(P_{0},P_{1},P_{1})=\left\{
R|f_{P_{0}P_{1}P_{1}}\left( R\right) =0\right\} ,  \label{b1.9}
\end{equation}
\begin{equation}
f_{P_{0}P_{1}P_{1}}\left( R\right) =S_{2}\left( P_{0},R,P_{1}\right) \equiv
\sqrt{2\sigma \left( P_{0},P_{1}\right) }-\sqrt{2\sigma \left(
P_{0},R\right) }-\sqrt{2\sigma \left( P_{1},R\right) }  \label{b1.10}
\end{equation}
In the proper Euclidean geometry ${\cal T}_{[P_{0}P_{1}]}$ is simply a
segment of the straight between the points $P_{0},P_{1}$.

The most important and interesting EGOs arise, when values of the envelope
function $f_{{\cal P}^n}(R)$ coincide with values of the function $F_{n+1}(%
{\cal P}^n,R)$, determined by relation (\ref{a1.6}) and proportional to the
squared length of the finite $\sigma$-subspace, consisting of $n+2$ points $%
{\cal P}^n,R$. This object is called the $n$th order natural geometric
object (NGO). It is defined by the relation (\ref{b1.3}). In the case of
proper Euclidean geometry it coincides with $n$-dimensional plane.

Another functions $f$ generate another envelopes of elementary geometrical
objects for the given skeleton ${\cal P}^n$.
For instance, the set of two points $\{P_0,P_1\}$ forms a skeleton not only
for the tube ${\cal T}_{P_0P_1}$, but also for the segment ${\cal T}%
_{[P_0P_1]}$ of the tube (straight) (\ref{b1.9}), and for the tube ray $%
{\cal T}_{[P_0P_1}$, which is defined by the relation
\begin{equation}
{\cal T}_{[P_0P_1}=\left\{ R|S_2\left( P_0,P_1,R\right) =0\right\}
\label{a2.9}
\end{equation}
where the function $S_2$ is defined by the relation (\ref{b1.10}).

\section{Interrelation between T-geometric and \\ Riemannian conceptions of
geometry}

\begin{opred}
The geometric conception is a totality of principles of the geometry
construction. \label{dd.1.1}
\end{opred}

Let us compare the Riemannian conception of geometry and that of T-geometry.
$n$-dimensional Riemannian geometry $R_n=\left\{ {\bf g},K,{\cal M}%
_n\right\} $ is introduced on $n$-dimensio\-nal manifold ${\cal M}_n$ in some
coordinate system $K$ by setting the metric tensor $g_{ik}(x)$, $%
i,k=1,2,\ldots n$. Thereafter, using the definition (\ref{b5.5}) of the
curve, which always can be introduced on the manifold ${\cal M}_n$, one
introduces concept of geodesic ${\cal L}_{[xx^{\prime }]}$ as the shortest
curve connecting points with coordinates $x$ and $x^{\prime }$. In the
Riemannian space $R_n=\left\{ {\bf g},K,{\cal M}_n\right\} $ one introduces
the world function $\sigma _R(x,x^{\prime })$ between points $x$ and $%
x^{\prime }$, defined by the relation
\begin{equation}
\sigma _R(x,x^{\prime })=\frac 12\left( \int\limits_{_{{\cal L}_{[xx^{\prime
}]}}}\sqrt{g_{ik}dx^idx^k}\right) ^2,  \label{a3.15}
\end{equation}
where ${\cal L}_{[xx^{\prime }]}$ denotes a segment of geodesic, connecting
points $x$ and $x^{\prime }$.

T-geometry can be introduced on any set $\Omega $, including the manifold $%
{\cal M}_{n}$. To set T-geometry on ${\cal M}_{n}$, it is insufficient of
the metric tensor $g_{ik}(x)$, $i,k=1,2,\ldots n$ introduction, because it
determines only first derivatives of world function at coinciding points
\begin{equation}
g_{ik}\left( x\right) =-\sigma _{ik^{\prime }}\left( x,x\right) \equiv -
\left[ \frac{\sigma \left( x,x^{\prime }\right) }{\partial x^{i}\partial
x^{\prime k}}\right] _{x^{\prime }=x}  \label{b6.1}
\end{equation}
This is insufficient for determination of the world function. For setting
T-geometry in a way consistent with the Riemannian geometry, one should set $%
\sigma (x,x^{\prime })=\sigma _{R}(x,x^{\prime })$, where $\sigma
_{R}(x,x^{\prime })$ is defined by the relation (\ref{a3.15}). Now one can
construct geometric objects by the method described above. The T-geometry,
introduced in such a way, will be referred to as $\sigma $-Riemannian
geometry, for distinguishing different conceptions (i.e. rules of
construction) of geometry.

Note that the world function, consistent with Riemannian geometry on the
manifold, may be set as a solution of equations in partial derivatives. For
instance, the world function can be defined as the solution of the
differential equation \cite{S60}
\begin{equation}
\sigma _{i}g^{ik}\left( x\right) \sigma _{k}=2\sigma ,\qquad \sigma
_{i}\equiv \frac{\partial \sigma }{\partial x^{i}} \qquad i=1,2,\ldots n,
\label{b6.2}
\end{equation}
satisfying the conditions (\ref{a2.1}).

The basic geometric objects of Riemannian geometry -- geodesic segments $%
{\cal L}_{[xx^{\prime }]}$ coincide with the first order NGOs in T-geometry
-- the tube segments ${\cal T}_{[xx^{\prime }]}$, defined by the relations (%
\ref{b1.9}). Thus one can say on partial coincidence of two geometric
conceptions: Riemannian and $\sigma $-Riemannian ones. But such a
coincidence is not complete. There are some difference which appears
sometimes essential.

Let us consider the case, when the manifold ${\cal M}_n$ coincides with $%
{\Bbb R}^n$ and metric tensor $g_{ik}=$const, $i,k=1,2,\ldots n$, $\;g=\det
||g_{ik}||\neq 0$ is the metric tensor of the proper Euclidean space. The
world function is described by the relation (\ref{a3.14}), and the proper
Riemannian space $E_n=\left\{ {\bf g}_E,K,{\Bbb R}^n\right\}$ is the proper
Euclidean space. Here ${\bf g}_E$ denotes the metric tensor of the proper
Euclidean space.

Now let us consider the proper Riemannian space $R_n=\left\{ {\bf g}%
_E,K,D\right\} $, where $D\subset {\Bbb R}^n$ is some region of the proper
Euclidean space $E_n=\left\{ {\bf g}_E,K,{\Bbb R}^n\right\} .$ If this
region $D$ is convex, i.e. any segment ${\cal L}_{[xx^{\prime }]}$ of
straight, passing through points $x,x^{\prime }\in D$, belongs to $D$ ($%
{\cal L}_{[xx^{\prime }]}\subset D)$, the world function of the proper
Riemannian space $R_n=\left\{ {\bf g}_E,K,D\right\} $ has the form (\ref
{a3.14}), and the proper Riemannian space $R_n=\left\{ {\bf g}_E,K,D\right\}
$ can be embedded isometrically to the proper Euclidean space $E_n=\left\{
{\bf g}_E,K,{\Bbb R}^n\right\} $.

If the region $D$ is not convex, the system of geodesics in the region $%
R_{n}=\left\{ {\bf g}_{E},K,D\right\} $ is not a system of straights, and
world function (\ref{a3.15}) is not described by the relation (\ref{a3.14}).
In this case the region $D$ cannot be embedded isometrically in $%
E_{n}=\left\{ {\bf g}_{E},K,{\Bbb R}^{n}\right\} $, in general. It seems to
be paradoxical that one (nonconvex) part of the proper Euclidean space
cannot be embedded isometrically to it, whereas another (convex) part can.

The convexity problem appears to be rather complicated, and most of
mathematicians prefer to go around this problem, dealing only with convex
regions \cite{A48}. In T-geometry there is no convexity problem. Indeed,
according to definition \ref{d3.1.1аe} subset of points of $\sigma $-space
is always embeddable isometrically in $\sigma $-space. From viewpoint of
T-geometry a removal of any region $R_{n}=\left\{ {\bf g}_{E},K,D\right\} $
from the proper Euclidean space $R_{n}=\left\{ {\bf g}_{E},K,{\Bbb R}%
^{n}\right\} $ cannot change shape of geodesics (first order NGOs). It leads
only to holes in geodesics, making them discontinuous. The continuity is a
property of the coordinate system, used in the proper Riemannian geometry as
a main tool of description. Using continuous coordinate systems for
description, we transfer constraints imposed on coordinate system to the
geometry itself.

Insisting on continuity of geodesics, one overestimates importance  of
continuity for geometry and attributes continuous geodesics (the first order
NGOs) to any proper Riemannian geometry, whereas the continuity is a special
property of the proper Euclidean geometry. From viewpoint of T-geometry the
convexity problem is an artificial problem. Existence of the convexity
problem in the Riemannian conception of geometry and its absence in
T-geometric conception means that the second conception of geometry is more
perfect.

\section{Riemannian geometry and one-dimensionality of the first order tubes}

Let us consider the $n$-dimensional pseudoeuclidean space $E_n=\left\{ {\bf %
g}_1,K,{\Bbb R}^n\right\} $ of the index $1$, ${\bf g}_1=$diag$\left\{
1,-1,-1\ldots -1\right\} $ to be a kind of $n$-dimensional Riemannian space.
The world function is defined by the relation (\ref{a3.14})
\begin{equation}
\sigma _1(x,x^{\prime })={\frac 12}\sum_{i,k=1}^ng^{ik}\left(
x_i-x_i^{\prime }\right) \left( x_k-x_k^{\prime }\right) ,\qquad g^{ik}={\rm %
diag}\left\{ 1,-1,-1\ldots -1\right\}  \label{c3.0}
\end{equation}
Geodesic ${\cal L}_{yy^{\prime }}$ is a straight line, and it is considered
in pseudoeuclidean geometry to be the first order NGOs, determined by two
points $y$ and $y^{\prime }$

\begin{equation}
{\cal L}_{yy^{\prime }}:\quad x^i=\left( y^i-y^{\prime i}\right) \tau
,\qquad i=1,2,\ldots n,\qquad \tau \in {\Bbb R}  \label{c3.1}
\end{equation}
The geodesic ${\cal L}_{yy^{\prime }}$ is called timelike, if $\sigma
_1(y,y^{\prime })>0$, and it is called spacelike if $\sigma _1(y,y^{\prime
})<0$. The geodesic ${\cal L}_{yy^{\prime }}$ is called null, if $\sigma
_1(y,y^{\prime })=0$.

The pseudoeuclidean space $E_n=\left\{ {\bf g}_1,K,{\Bbb R}^n\right\} $
generates the $\sigma $-space $V=\left\{ \sigma _1,{\Bbb R}^n\right\} $,
where the world function $\sigma _1$ is defined by the relation (\ref{c3.0}%
). The first order tube (NGO) ${\cal T}\left( x,x^{\prime }\right) $ in the $%
\sigma $-Riemannian space $V=\left\{ \sigma _1,{\Bbb R}^n\right\} $ is
defined by the relation (\ref{b1.3})
\begin{equation}
{\cal T}\left( x,x^{\prime }\right) \equiv {\cal T}_{xx^{\prime }}=\left\{
r|F_2\left( x,x^{\prime },r\right) =0\right\} ,\qquad \sigma _1(x,x^{\prime
})\neq 0,\qquad x,x^{\prime },r\in {\Bbb R}^n,  \label{c3.2}
\end{equation}
\begin{equation}
F_2\left( x,x^{\prime },r\right) =\left|
\begin{array}{cc}
(x_i^{\prime }-x_i)(x^{\prime i}-x^i) & (x_i^{\prime }-x_i)(r^i-x^i) \\
(r_i-x_i)(x^{\prime i}-x^i) & (r_i-x_i)(r^i-x^i)
\end{array}
\right|  \label{c3.3}
\end{equation}
Solution of equations (\ref{c3.2}), (\ref{c3.3}) gives the following result
\begin{equation}
{\cal T}_{xx^{\prime }}=\left\{ r\left| \bigcup\limits_{y\in {\Bbb R}%
^n}\bigcup\limits_{\tau \in {\Bbb R}}r=\left( x^{\prime }-x\right) \tau
+y-x\wedge \Gamma (x,x^{\prime },y)=0\right. \wedge \Gamma (x,y,y)=0\right\}
,  \label{c3.4}
\end{equation}
\[
x,x^{\prime },y,r\in {\Bbb R}^n
\]
where $\Gamma (x,x^{\prime },y)=(x_i^{\prime }-x_i)(y^i-x^i)$ is the scalar
product of vectors $\overrightarrow{xy}$ and $\overrightarrow{xx^{\prime }}$
defined by the relation (\ref{a1.7}). In the case of timelike vector $%
\overrightarrow{xx^{\prime }}$, when $\sigma _1(x,x^{\prime })>0$, there is
a unique null vector $\overrightarrow{xy}=\overrightarrow{xx}=%
\overrightarrow{0}$ which is orthogonal to the vector $\overrightarrow{%
xx^{\prime }}$. In this case the ($n-1)$-dimensional surface ${\cal T}%
_{xx^{\prime }}$ degenerates into the one-dimensional straight
\begin{equation}
{\cal T}_{xx^{\prime }}=\left\{ r\left| \bigcup\limits_{\tau \in {\Bbb R}%
}r=\left( x^{\prime }-x\right) \tau \right. \right\} ,\qquad \sigma
_1(x,x^{\prime })>0,\qquad x,x^{\prime },r\in {\Bbb R}^n,  \label{c3.5}
\end{equation}
Thus, for timelike vector $\overrightarrow{xx^{\prime }}$ the first order
tube ${\cal T}_{xx^{\prime }}$ coincides with the geodesic ${\cal L}%
_{xx^{\prime }}$. In the case of spacelike vector $\overrightarrow{%
xx^{\prime }}$ the $(n-1)$-dimensional tube ${\cal T}_{xx^{\prime }}$
contains the one-dimensional geodesic ${\cal L}_{xx^{\prime }}$ of the
pseudoeuclidean space $E_n=\left\{ {\bf g}_1,K,{\Bbb R}^n\right\} $.

This difference poses the question what is the reason of this difference and
what of the two generalization of the proper Euclidean geometry is more
reasonable. Note that four-dimensional pseudoeuclidean geometry is used for
description of the real space-time. One can try to resolve this problem from
experimental viewpoint. Free classical particles are described by means of
timelike straight lines. At this point the pseudoeuclidean geometry and the
$\sigma $-pseudoeuclidean geometry (T-geometry) lead to the same result.
The spacelike straights are believed to describe the particles moving with
superlight speed (so-called taxyons). Experimental attempts of taxyons
discovery were failed. Of course, trying to discover taxyons, one considered
them to be described by spacelike straights. On the other hand, the
physicists believe that all what can exist does exist and may be discovered.
From this viewpoint the failure of discovery of taxyons in the form of
spacelike line justifies in favor of taxyons in the form of
three-dimensional surfaces.

To interpret the structure of the set (\ref{c3.4}), describing the first
order tube, let us take into account the zeroth order tube ${\cal T}_x$,
determined by the point $x$ in the $\sigma $-pseudoeuclidean space is the
light cone with the vertex at the point $x$ (not the point $x$). Practically
the first order tube consists of such sections of the light cones with their
vertex $y\in {\cal L}_{xx^{\prime }}$ that all vectors $\overrightarrow{yr}$
of these sections are orthogonal to the vector $\overrightarrow{xx^{\prime }}
$. In other words, the first order tube ${\cal T}_{xx^{\prime }}$ consists
of the zeroth order tubes ${\cal T}_y$ sections at $y$, orthogonal to $%
\overrightarrow{xx^{\prime }}$, with $y\in {\cal L}_{xx^{\prime }}$. For
timelike $\overrightarrow{xx^{\prime }}$ this section consists of one point,
but for the spacelike $\overrightarrow{xx^{\prime }}$ it is two-dimensional
section of the light cone.

\section{Collinearity in Riemannian and $\protect\sigma $-Riemannian geometry%
}

Let us return to the Riemannian space $R_n=\left\{ {\bf g},K,D\right\}
,\,\quad D\subset {\Bbb R}^n$, which generates the world function $\sigma
(x,x^{\prime })$ defined by the relation (\ref{a3.15}). Then the $\sigma $%
-space $V=\left\{ \sigma ,D\right\} $ appears. it will be referred to as $%
\sigma $-Riemannian space. We are going to compare concept of collinearity
(parallelism) of two vectors in the two spaces.

The world function $\sigma =\sigma (x,x^{\prime })$ of both $\sigma $%
-Riemannian and Riemannian spaces satisfies the system of equations \cite
{R62}\footnote{%
The paper \cite{R62} is hardly available for English speaking reader. Survey
of main results of \cite{R62} in English may be found in \cite{R64}. See
also \cite{R92}}
\begin{equation}
\begin{array}{cc}
(1)\quad \sigma _l\sigma ^{{lj^{\prime }}}\sigma _{j^{\prime }}=2\sigma
\quad \quad \quad & (4)\quad \det \parallel \sigma _{i||k}\parallel \neq 0
\\
(2)\quad \sigma (x,x^{\prime })=\sigma (x^{\prime },x)\quad \quad & (5)\quad
\det \parallel \sigma _{{ik^{\prime }}}\parallel \neq 0 \\
(3)\quad \sigma (x,x)=0\quad \quad \quad \quad \quad & (6)\quad \sigma
_{i||k||l}=0\quad \quad
\end{array}
\label{b2.39}
\end{equation}
where the following designations are used
\[
\sigma _i\equiv \frac{\partial \sigma }{\partial x^i},\qquad \sigma
_{i^{\prime }}\equiv \frac{\partial \sigma }{\partial x^{\prime i}},\qquad
\sigma _{ik^{\prime }}\equiv \frac{\partial ^2\sigma }{\partial x^i\partial
x^{\prime k}},\qquad \sigma ^{ik^{\prime }}\sigma _{lk^{\prime }}=\delta
_l^i
\]
Here the primed index corresponds to the point $x'$, and unprimed index
corresponds to the point $x$. Two parallel vertical strokes mean covariant
derivative $\tilde \nabla _i^{x^{\prime }}$ with respect to $x^i$ with the
Christoffel symbol
\[
\Gamma _{kl}^i\equiv \Gamma _{kl}^i\left( x,x^{\prime }\right) \equiv \sigma
^{is^{\prime }}\sigma _{kls^{\prime }},\qquad \sigma _{kls^{\prime }}\equiv
\frac{\partial ^3\sigma }{\partial x^k\partial x^l\partial x^{\prime s}}
\]
For instance,
\begin{equation}
G_{ik}\equiv G_{ik}(x,x^{\prime })\equiv \sigma _{i||k}\equiv \frac{\partial
\sigma _i}{\partial x^k}-\Gamma _{ik}^l\left( x,x^{\prime }\right) \sigma
_l\equiv \frac{\partial \sigma _i}{\partial x^k}-\sigma _{iks^{\prime
}}\sigma ^{ls^{\prime }}\sigma _l  \label{a2.0}
\end{equation}
\[
G_{ik||l}\equiv \frac{\partial G_{ik}}{\partial x^l}-\sigma _{ils^{\prime
}}\sigma ^{js^{\prime }}G_{jk}-\sigma _{kls^{\prime }}\sigma ^{js^{\prime
}}G_{ij}
\]
Summation from $1$ to $n$ is produced over repeated indices. The covariant
derivative $\tilde \nabla _i^{x^{\prime }}$ with respect to $x^i$ with the
Christoffel symbol $\Gamma _{kl}^i\left( x,x^{\prime }\right) $ acts only on
the point $x$ and on unprimed indices. It is called the tangent derivative,
because it is a covariant derivative in the Euclidean space $E_{x^{\prime }}$
which is tangent to the Riemannian space $R_n$ at the point $x^{\prime }$.
The covariant derivative $\tilde \nabla _{i^{\prime }}^x$ with respect to $%
x^{\prime i}$ with the Christoffel symbol $\Gamma _{k^{\prime }l^{\prime
}}^{i^{\prime }}\left( x,x^{\prime }\right) $ acts only on the point $%
x^{\prime }$ and on primed indices. It is a covariant derivative in the
Euclidean space $E_x$ which is tangent to the $\sigma $-Riemannian space $%
R_n $ at the point $x$ \cite{R62}.

In general, the world function $\sigma $ carries out the geodesic mapping $%
G_{x^{\prime }}:$ $R_n\rightarrow E_{x^{\prime }}$ of the Riemannian space $%
R_n=\left\{ {\bf g},K,D\right\} $ on the Euclidean space $E_{x^{\prime
}}=\left\{ {\bf g},K_{x^{\prime }},D\right\} $, tangent to $R_n=\left\{ {\bf %
g},K,D\right\} $ at the point $x^{\prime }$ \cite{R62}. This mapping
transforms the coordinate system $K$ in $R_n$ into the coordinate system $%
K_{x^{\prime }}$ in $E_{x^{\prime }}$. The mapping is geodesic in the sense
that it conserves the lengths of segments of all geodesics, passing through
the tangent point $x^{\prime }$ and angles between them at this point.

The tensor $G_{ik}$, defined by (\ref{a2.0}) is the metric tensor at the
point $x$ in the tangent Euclidean space $E_{x^{\prime }}$. The covariant
derivatives $\tilde \nabla _i^{x^{\prime }}$ and $\tilde \nabla
_k^{x^{\prime }}$ commute identically, i.e. $(\tilde \nabla _i^{x^{\prime }}
\tilde \nabla _k^{x^{\prime }}-\tilde \nabla _k^{x^{\prime }}\tilde \nabla
_i^{x^{\prime }})A_{ls}\equiv 0,$ for any tensor $A_{ls}$ \cite{R62}. This
shows that they are covariant derivatives in the flat space $E_{x^{\prime }}$%
.

The system of equations (\ref{b2.39}) contains only world function $\sigma $
and its derivatives, nevertheless the system of equations (\ref{b2.39}) is
not $\sigma $-immanent, because it contains a reference to a coordinate
system. It does not contain the metric tensor explicitly. Hence, it is valid
for any Riemannian space $R_n=\left\{ {\bf g},K,D\right\} $. All relations
written above are valid also for the $\sigma $-space $V=\left\{ \sigma
,D\right\} $, provided the world function $\sigma $ is coupled with the
metric tensor by relation (\ref{a3.15}).

$\sigma $-immanent expression for scalar product $\left( {\bf P}_0{\bf P}_1.%
{\bf Q}_0{\bf Q}_1\right) $ of two vectors ${\bf P}_0{\bf P}_1$ and ${\bf Q}%
_0{\bf Q}_1$ in the proper Euclidean space has the form
\begin{equation}
\left( {\bf P}_0{\bf P}_1.{\bf Q}_0{\bf Q}_1\right) \equiv \sigma \left(
P_0,Q_1\right) +\sigma \left( Q_0,P_1\right) -\sigma \left( P_0,Q_0\right)
-\sigma \left( P_1,Q_1\right)  \label{a2.1a}
\end{equation}
This relation can be easily proved as follows.

In the proper Euclidean space three vectors ${\bf P}_0{\bf P}_1$, ${\bf P}_0%
{\bf Q}_1$, and ${\bf P}_1{\bf Q}_1$ are coupled by the relation
\begin{equation}
\mid {\bf P}_1{\bf Q}_1\mid ^2=\mid {\bf P}_0{\bf Q}_1-{\bf P}_0{\bf P}%
_1\mid ^2=\mid {\bf P}_0{\bf P}_1\mid ^2+\mid {\bf P}_0{\bf Q}_1\mid ^2-2(%
{\bf P}_0{\bf P}_1.{\bf P}_0{\bf Q}_1)  \label{a2.15}
\end{equation}
where $({\bf P}_0{\bf P}_1.{\bf P}_0{\bf Q}_1)$ denotes the scalar product
of two vectors ${\bf P}_0{\bf P}_1$ and ${\bf P}_0{\bf Q}_1$ in the proper
Euclidean space. It follows from (\ref{a2.15})
\begin{equation}
({\bf P}_0{\bf P}_1.{\bf P}_0{\bf Q}_1)={\frac 12}\{\mid {\bf P}_0{\bf Q}%
_1\mid ^2+\mid {\bf P}_0{\bf P}_1\mid ^2-\mid {\bf P}_1{\bf Q}_1\mid ^2\}
\label{a2.16}
\end{equation}
Substituting the point $Q_1$ by $Q_0$ in (\ref{a2.16}), one obtains
\begin{equation}
({\bf P}_0{\bf P}_1.{\bf P}_0{\bf Q}_0)={\frac 12}\{\mid {\bf P}_0{\bf Q}%
_0\mid ^2+\mid {\bf P}_0{\bf P}_1\mid ^2-\mid {\bf P}_1{\bf Q}_0\mid ^2\}
\label{a2.16a}
\end{equation}
Subtracting (\ref{a2.16a}) from (\ref{a2.16}) and using the properties of
the scalar product in the proper Euclidean space, one obtains
\begin{equation}
({\bf P}_0{\bf P}_1.{\bf Q}_0{\bf Q}_1)={\frac 12}\{\mid {\bf P}_0{\bf Q}%
_1\mid ^2+\mid {\bf Q}_0{\bf P}_1\mid ^2-\mid {\bf P}_0{\bf Q}_0\mid ^2-\mid
{\bf P}_1{\bf Q}_1\mid ^2\}  \label{a2.16b}
\end{equation}
Taking into account that $\mid {\bf P}_0{\bf Q}_1\mid ^2=2\sigma \left(
P_0,Q_1\right) $, one obtains the relation\ (\ref{a2.1a}) from the relation (%
\ref{a2.16b}).

Two vectors ${\bf P}_0{\bf P}_1$ and ${\bf Q}_0{\bf Q}_1$ are collinear $%
{\bf P}_0{\bf P}_1||{\bf Q}_0{\bf Q}_1$ (parallel or antiparallel), provided
$\cos ^2\theta =1,$ where $\theta $ is the angle between the vectors ${\bf P}%
_0{\bf P}_1$ and ${\bf Q}_0{\bf Q}_1$. Taking into account that

\begin{equation}
\cos ^2\theta =\frac{\left( {\bf P}_0{\bf P}_1.{\bf Q}_0{\bf Q}_1\right) ^2}{%
\left( {\bf P}_0{\bf P}_1.{\bf P}_0{\bf P}_1\right) ({\bf Q}_0{\bf Q}_1.{\bf %
Q}_0{\bf Q}_1)}=\frac{\left( {\bf P}_0{\bf P}_1.{\bf Q}_0{\bf Q}_1\right) ^2%
}{|{\bf P}_0{\bf P}_1|^2\cdot |{\bf Q}_0{\bf Q}_1|^2}  \label{a2.2}
\end{equation}
one obtains the following $\sigma $-immanent condition of the two vectors
collinearity
\begin{equation}
{\bf P}_0{\bf P}_1||{\bf Q}_0{\bf Q}_1:\qquad \left( {\bf P}_0{\bf P}_1.{\bf %
Q}_0{\bf Q}_1\right) ^2=|{\bf P}_0{\bf P}_1|^2\cdot |{\bf Q}_0{\bf Q}_1|^2
\label{a2.3}
\end{equation}
The collinearity condition (\ref{a2.3}) is $\sigma $-immanent, because by
means of (\ref{a2.1a}) it can be written in terms of the $\sigma $-function
only. Thus, this relation describes the vectors collinearity in the case of
arbitrary $\sigma $-space.

Let us describe this relation for the case of $\sigma $-Riemannian geometry.
Let coordinates of the points $P_0,P_1,Q_0,Q_1$ be respectively $x,$ $x+dx,$
$x^{\prime }$ and $x^{\prime }+dx^{\prime }$. Then writing (\ref{a2.1a}) and
expanding it over $dx$ and $dx^{\prime }$, one obtains
\begin{eqnarray*}
\left( {\bf P}_0{\bf P}_1.{\bf Q}_0{\bf Q}_1\right) &\equiv &\sigma \left(
x,x^{\prime }+dx^{\prime }\right) +\sigma \left( x^{\prime },x+dx\right)
-\sigma \left( x,x^{\prime }\right) -\sigma \left( x+dx,x^{\prime
}+dx^{\prime }\right)  \label{a2.1b} \\
& = & \sigma _{l^{\prime }}dx^{\prime l^{\prime }}+\frac 12\sigma
_{l^{\prime },s^{\prime }}dx^{\prime l^{\prime }}dx^{\prime s^{\prime }}
+\sigma _idx^i+\frac 12\sigma _{i,k}dx^idx^k \\
&&\ -\sigma _idx^i-\sigma _{l^{\prime }}dx^{\prime l^{\prime }}-\frac
12\sigma _{i,k}dx^idx^k-\sigma _{i,l^{\prime }}dx^idx^{\prime l^{\prime
}}-\frac 12\sigma _{l^{\prime },s^{\prime }}dx^{\prime l^{\prime
}}dx^{\prime s^{\prime }}
\end{eqnarray*}
\begin{equation}
\left( {\bf P}_0{\bf P}_1.{\bf Q}_0{\bf Q}_1\right) =-\sigma _{i,l^{\prime
}}dx^idx^{\prime l^{\prime }}=-\sigma _{il^{\prime }}dx^idx^{\prime
l^{\prime }}  \label{a2.4}
\end{equation}
Here comma means differentiation. For instance, $\sigma _{i,k}\equiv
\partial \sigma _i/\partial x^k$. One obtains for $|{\bf P}_0{\bf P}_1|^2$
and $|{\bf Q}_0{\bf Q}_1|^2$%
\begin{equation}
|{\bf P}_0{\bf P}_1|^2=g_{ik}dx^idx^k,\qquad |{\bf Q}_0{\bf Q}%
_1|^2=g_{l^{\prime }s^{\prime }}dx^{\prime l^{\prime }}dx^{\prime s^{\prime
}}  \label{a2.5}
\end{equation}
where $g_{ik}=g_{ik}(x)$ and $g_{l^{\prime }s^{\prime }}=g_{l^{\prime
}s^{\prime }}(x^{\prime })$. Then the collinearity condition (\ref{a2.3}) is
written in the form
\begin{equation}
\left( \sigma _{il^{\prime }}\sigma _{ks^{\prime }}-g_{ik}g_{l^{\prime
}s^{\prime }}\right) dx^idx^kdx^{\prime l^{\prime }}dx^{\prime s^{\prime }}=0
\label{a4.3}
\end{equation}
Let us take into account that in the Riemannian space the metric tensor $%
g_{l^{\prime }s^{\prime }}$ at the point $x^{\prime }$ can be expressed via
the world function $\sigma $ of points $x,x^{\prime }$ by means of the
relation \cite{R62}
\begin{equation}
g_{l^{\prime }s^{\prime }}=\sigma _{il^{\prime }}G^{ik}\sigma _{ks^{\prime
}},\qquad g^{l^{\prime }s^{\prime }}=\sigma ^{il^{\prime }}G_{ik}\sigma
^{ks^{\prime }}  \label{a4.4}
\end{equation}
where the tensor $G_{ik}$ is defined by the relation (\ref{a2.0}), and $%
G^{ik}$ is defined by the relation
\begin{equation}
G^{il}G_{lk}=\delta _k^i  \label{a4.5}
\end{equation}

Substituting the first relation (\ref{a2.0}) in (\ref{a4.3}) and using
designation
\begin{equation}
u_{i}=-\sigma _{il^{\prime }}dx^{\prime l^{\prime }},\qquad
u^{i}=G^{ik}u_{k}=-\sigma ^{il^{\prime }}g_{l^{\prime }s^{\prime
}}dx^{s^{\prime }}  \label{a4.6}
\end{equation}
one obtains
\begin{equation}
\left( \delta _{i}^{l}\delta _{k}^{s}-g_{ik}G^{ls}\right)
u_{l}u_{s}dx^{i}dx^{k}=0  \label{a2.6}
\end{equation}
The vector $u_{i}$ is the vector $dx_{i^{\prime }}^{\prime }=g_{i^{\prime
}k^{\prime }}dx^{\prime k^{\prime }}$ transported parallelly from the point $%
x^{\prime }$ to the point $x$ in the Euclidean space $E_{x^{\prime }}$
tangent to the Riemannian space $R_{n}$. Indeed,
\begin{equation}
u_{i}=-\sigma _{il^{\prime }}g^{l^{\prime }s^{\prime }}dx_{s}^{\prime
},\qquad \tilde{\nabla}_{k}^{x^{\prime }}\left( -\sigma _{il^{\prime
}}g^{l^{\prime }s^{\prime }}\right) \equiv 0,\qquad i,k=1,2,\ldots n
\label{a2.7}
\end{equation}
and tensor $-\sigma _{il^{\prime }}g^{l^{\prime }s^{\prime }}$ is the
operator of the parallel transport in $E_{x^{\prime }}$, because
\[
\left[ -\sigma _{il^{\prime }}g^{l^{\prime }s^{\prime }}\right]
_{x=x^{\prime }}=\delta _{i^{\prime }}^{s^{\prime }}
\]
and the tangent derivative of this operator is equal to zero identically.
For the same reason, i.e. because of
\[
\left[ \sigma ^{il^{\prime }}g_{l^{\prime }s^{\prime }}\sigma ^{ks^{\prime }}%
\right] _{x=x^{\prime }}=g^{i^{\prime }k^{\prime }},\qquad \tilde{\nabla}%
_{s}^{x^{\prime }}(\sigma ^{il^{\prime }}g_{l^{\prime }s^{\prime }}\sigma
^{ks^{\prime }})\equiv 0
\]
$G^{ik}=\sigma ^{il^{\prime }}g_{l^{\prime }s^{\prime }}\sigma ^{ks^{\prime
}}$ is the contravariant metric tensor in $E_{x^{\prime }}$, at the point $x$%
.

The relation (\ref{a2.6}) contains vectors at the point $x$ only . At fixed $%
u_i=-\sigma _{il^{\prime }}dx^{\prime l^{\prime }}$ it describes a
collinearity cone, i.e. a cone of infinitesimal vectors $dx^i$ at the point $%
x $ parallel to the vector $dx^{\prime i^{\prime }}$ at the point $x^{\prime
}$. Under some condition the collinearity cone can degenerates into a line.
In this case there is only one direction, parallel to the fixed vector $u^i$%
. Let us investigate, when this situation takes place.

At the point $x$ two metric tensors $g_{ik}$ and $G_{ik}$ are connected by
the relation \cite{R62}
\begin{equation}
G_{ik}(x,x^{\prime })=g_{ik}(x)+\int\limits_x^{x^{\prime }}F_{{ikj^{\prime
\prime }}s^{{\prime \prime }}}(x,x^{\prime \prime })\sigma ^{j^{{\prime
\prime }}}(x,x^{\prime \prime })d{x^{\prime \prime }}^{s^{\prime \prime }},
\label{b3.23}
\end{equation}
where according to \cite{R62}
\begin{equation}
\sigma ^{i^{\prime }}=\sigma ^{li^{\prime }}\sigma _l=G^{l^{\prime
}i^{\prime }}\sigma _{l^{\prime }}=g^{l^{\prime }i^{\prime }}\sigma
_{l^{\prime }}  \label{a4.6a}
\end{equation}
Integration does not depend on the path, because it is produced in the
Euclidean space $E_{x^{\prime }}$. The two-point tensor $F_{{ilk^{\prime
}j^{\prime }}}=F_{{ilk^{\prime }j^{\prime }}}(x,x^{\prime })$ is the
two-point curvature tensor, defined by the relation
\begin{equation}
F_{{ilk^{\prime }j^{\prime }}}=\sigma _{{ilj^{\prime }}\parallel k^{\prime
}}=\sigma _{{ilj^{\prime }},k^{\prime }}-\sigma _{{sj^{\prime }k^{\prime }}%
}\sigma ^{{sm^{\prime }}}\sigma _{{ilm^{\prime }}}=\sigma _{i\mid
l||k^{\prime }||j^{\prime }}  \label{b3.5}
\end{equation}
where one vertical stroke denotes usual covariant derivative and two
vertical strokes denote tangent derivative. The two-point curvature tensor $%
F_{{ilk^{\prime }j^{\prime }}}$ has the following symmetry properties
\begin{equation}
F_{{ilk^{\prime }j^{\prime }}}=F_{{lik^{\prime }j^{\prime }}}=F_{{%
ilj^{\prime }k^{\prime }}},\qquad F_{{ilk^{\prime }j^{\prime }}}(x,x^{\prime
})=F_{{k^{\prime }j^{\prime }il}}(x^{\prime },x)  \label{b3.10}
\end{equation}
It is connected with the one-point Riemann-Ghristoffel curvature tensor $%
r_{iljk}$ by means of relations
\begin{equation}
r_{{iljk}}=\left[ F_{{ikj^{\prime }l^{\prime }}}-F_{{ijk^{\prime }l^{\prime }%
}}\right] _{x^{\prime }=x}=f_{{ikjl}}-f_{{ijkl}},\qquad f_{{iklj}}=\left[ F_{%
{ikj^{\prime }l^{\prime }}}\right] _{x^{\prime }=x}  \label{b3.14}
\end{equation}

In the Euclidean space the two-point curvature tensor $F_{{ilk^{\prime
}j^{\prime }}}$ vanishes as well as the Riemann-Ghristoffel curvature tensor
$r_{iljk}$.

Let us introduce designation
\begin{equation}
\Delta _{ik}=\Delta _{ik}(x,x^{\prime })=\int\limits_x^{x^{\prime }}F_{{%
ikj^{\prime \prime }}s^{{\prime \prime }}}(x,x^{\prime \prime })\sigma ^{j^{{%
\prime \prime }}}(x,x^{\prime \prime })d{x^{\prime \prime }}^{s^{\prime
\prime }}  \label{a4.8}
\end{equation}
and choose the geodesic ${\cal L}_{xx^{\prime }}$ as the path of
integration. It is described by the relation
\begin{equation}
\sigma _i(x,x^{\prime \prime })=\tau \sigma _i(x,x^{\prime })  \label{a4.11}
\end{equation}
which determines $x^{\prime \prime }$ as a function of parameter $\tau $.
Differentiating with respect to $\tau $, one obtains
\begin{equation}
\sigma _{ik^{\prime \prime }}(x,x^{\prime \prime })dx^{\prime \prime
k^{\prime \prime }}=\sigma _i(x,x^{\prime })d\tau  \label{a4.10}
\end{equation}
Resolving equations (\ref{a4.10}) with respect to $dx^{\prime \prime }$ and
substituting in (\ref{a4.8}), one obtains
\begin{equation}
\Delta _{ik}(x,x^{\prime })=\sigma _l(x,x^{\prime })\sigma _p(x,x^{\prime
})\int\limits_0^1F_{{ikj^{\prime \prime }}s^{{\prime \prime }}}(x,x^{\prime
\prime })\sigma ^{lj^{{\prime \prime }}}(x,x^{\prime \prime })\sigma ^{ps^{{%
\prime \prime }}}(x,x^{\prime \prime })\tau d{\tau }  \label{a4.12}
\end{equation}
where $x^{\prime \prime }$ is determined from (\ref{a4.11}) as a function of
$\tau $. Let us set
\begin{equation}
F_{{ik}}^{\;..lp}(x,x^{\prime })=F_{{ikj^{\prime }}s^{{\prime }%
}}(x,x^{\prime })\sigma ^{lj^{{\prime }}}(x,x^{\prime })\sigma ^{ps^{{\prime
}}}(x,x^{\prime })  \label{a4.13}
\end{equation}
then
\begin{equation}
G_{ik}(x,x^{\prime })=g_{ik}(x)+\Delta _{ik}(x,x^{\prime })  \label{a4.14}
\end{equation}
\begin{equation}
\Delta _{ik}(x,x^{\prime })=\sigma _l(x,x^{\prime })\sigma _p(x,x^{\prime
})\int\limits_0^1F_{{ik}}^{\;..lp}(x,x^{\prime \prime })\tau d{\tau }
\label{a4.15}
\end{equation}
Substituting $g_{ik}$ from (\ref{a4.14}) in (\ref{a2.6}), one obtains
\begin{equation}
\left( \delta _i^l\delta _k^s-G^{ls}\left( G_{ik}-\Delta _{ik}\right)
\right) u_lu_sdx^idx^k=0  \label{a4.7}
\end{equation}

Let us look for solutions of equation in the form of expansion
\begin{equation}
dx^i=\alpha u^i+v^i,\qquad G_{ik}u^iv^k=0  \label{a4.16}
\end{equation}
Substituting (\ref{a4.16}) in (\ref{a4.7}), one obtains equation for $v^i$%
\begin{equation}
G_{ls}u^lu^s\left[ G_{ik}v^iv^k-\Delta _{ik}\left( \alpha u^i+v^i\right)
\left( \alpha u^k+v^k\right) \right] =0  \label{a4.18}
\end{equation}
If the $\sigma $-Riemannian space $V=\left\{ \sigma ,D\right\} $ is $\sigma $%
-Euclidean, then as it follows from (\ref{a4.15}) $\Delta _{ik}=0$. If $%
V=\left\{ \sigma ,D\right\} $ is the proper $\sigma $-Euclidean space, $%
G_{ls}u^lu^s\neq 0$, and one obtains two equations for determination of $v^i$
\begin{equation}
G_{ik}v^iv^k=0,\qquad G_{ik}u^iv^k=0  \label{a4.17}
\end{equation}
The only solution
\begin{equation}
v^i=0,\qquad dx^i=\alpha u^i,\qquad i=1,2,\ldots n  \label{a2.6b}
\end{equation}
of (\ref{a4.18}) is a solution of the equation (\ref{a4.7}), where $\alpha $
is an arbitrary constant. In the proper Euclidean geometry the collinearity
cone always degenerates into a line.

Let now the space $V=\left\{ \sigma ,D\right\} $ be the $\sigma $%
-pseudoeuclidean space of index $1$, and the vector $u^i$ be timelike, i.e.
$G_{ik}u^iu^k>0$. Then equations (\ref{a4.17}) also have the solution (\ref
{a2.6b}). If the vector $u^i$ is spacelike, $G_{ik}u^iu^k<0,$ then two
equations (\ref{a4.17}) have non-trivial solution, and the collinearity cone
does not degenerate into a line. The collinearity cone is a section of the
light cone $G_{ik}v^iv^k=0$ by the plane $G_{ik}u^iv^k=0$. If the vector $%
u^i $ is null, $G_{ik}u^iu^k=0,$ then equation (\ref{a4.18}) reduces to the
form
\begin{equation}
G_{ik}u^iu^k=0,\qquad G_{ik}u^iv^k=0  \label{a4.19}
\end{equation}
In this case (\ref{a2.6b}) is a solution, but besides there are spacelike
vectors $v^i$ which are orthogonal to null vector $u^i$ and the collinearity
cone does not degenerate into a line.

In the case of the proper $\sigma $-Riemannian space $G_{ik}u^iu^k>0,$ and
equation (\ref{a4.18}) reduces to the form
\begin{equation}
G_{ik}v^iv^k-\Delta _{ik}\left( \alpha u^i+v^i\right) \left( \alpha
u^k+v^k\right) =0  \label{a4.20}
\end{equation}
In this case $\Delta _{ik}\neq 0$ in general, and the collinearity cone does
not degenerate. $\Delta _{ik}$ depends on the curvature an on the distance
between the points $x$ and $x^{\prime }$. The more space curvature and the
distance $\rho (x,x^{\prime }),$ the more the collinearity cone aperture.

In the curved proper $\sigma $-Riemannian space there is an interesting
special case, when the collinearity cone degenerates . In any $\sigma $%
-Riemannian space the following equality takes place \cite{R62}
\begin{equation}
G_{ik}\sigma ^k=g_{ik}\sigma ^k,\qquad \sigma ^k\equiv g^{kl}\sigma _l
\label{a4.21}
\end{equation}
Then it follows from (\ref{a4.14}) that
\begin{equation}
\Delta _{ik}\sigma ^k=0  \label{a4.22}
\end{equation}
It means that in the case, when the vector $u^i$ is directed along the
geodesic, connecting points $x$ and $x^{\prime }$, i.e. $u^i=\beta \sigma ^i$%
, the equation (\ref{a4.20}) reduces to the form
\begin{equation}
\left( G_{ik}-\Delta _{ik}\right) v^iv^k=0,\qquad u^i=\beta \sigma ^i
\label{a4.23}
\end{equation}
If $\Delta _{ik}$ is small enough as compared with $G_{ik},$ then
eigenvalues of the matrix $G_{ik}-\Delta _{ik}$ have the same sign, as those
of the matrix $G_{ik}.$ In this case equation (\ref{a4.23}) has the only
solution (\ref{a2.6b}), and the collinearity cone degenerates.

\section{Discussion}

Thus, we see that in the $\sigma $-Riemannian geometry at the point $x$
there are many vectors parallel to given vector at the point $x^{\prime }$.
This set of parallel vectors is described by the collinearity cone.
Degeneration of the collinearity cone into a line, when there is only one
direction, parallel to the given direction, is an exception rather than a
rule, although in the proper Euclidean geometry this degeneration takes
place always. Nonuniformity of space destroys the collinearity cone
degeneration. In the proper Riemannian geometry, where the world function
satisfies the system (\ref{b2.39}), one succeeded in conserving this
degeneration for direction along the geodesic, connecting points $x$ and $%
x^{\prime }$. This circumstance is very important for degeneration of the
first order NGOs into geodesic, because degeneration of NGOs is connected
closely with the collinearity cone degeneration.

Indeed, definition of the first order tube (\ref{b1.3}), or (\ref{c3.2}) may
be written also in the form
\begin{equation}
{\cal T}\left( {\cal P}^1\right) \equiv {\cal T}_{P_0P_1}=\left\{ R\left| \;
{\bf P}_0{\bf P}_1||{\bf P}_0{\bf R}\right. \right\} ,\qquad P_0,P_1,R\in
\Omega ,  \label{b5.1}
\end{equation}
where collinearity ${\bf P}_0{\bf P}_1||{\bf P}_0{\bf R}$ of two vectors $%
{\bf P}_0{\bf P}_1$ and ${\bf P}_0{\bf R}$ is defined by the $\sigma $%
-immanent relation (\ref{a2.3}), which can be written in the form

\begin{equation}
{\bf P}_0{\bf P}_1||{\bf P}_0{\bf R}:\quad F_2\left( P_0,P_1,R\right)
=\left|
\begin{array}{cc}
\left( {\bf P}_0{\bf P}_1.{\bf P}_0{\bf P}_1\right) & \left( {\bf P}_0{\bf P}%
_1.{\bf P}_0{\bf R}\right) \\
\left( {\bf P}_0{\bf R}.{\bf P}_0{\bf P}_1\right) & \left( {\bf P}_0{\bf R}.%
{\bf P}_0{\bf R}\right)
\end{array}
\right| =0  \label{b5.2}
\end{equation}
The form (\ref{b5.1}) of the first order tube definition allows one to
define the first order tube ${\cal T}(P_0,P_1;Q_0)$, passing through the
point $Q_0$ collinear to the given vector ${\bf P}_0{\bf P}_1$. This
definition has the $\sigma $-immanent form
\begin{equation}
{\cal T}(P_0,P_1;Q_0)=\left\{ R\left|\; {\bf P}_0{\bf P}_1||{\bf Q}_0{\bf R}%
\right. \right\} ,\qquad P_0,P_1,Q_0,R\in \Omega ,  \label{b5.3}
\end{equation}
where collinearity ${\bf P}_0{\bf P}_1||{\bf Q}_0{\bf R}$ of two vectors $%
{\bf P}_0{\bf P}_1$ and ${\bf Q}_0{\bf R}$ is defined by the $\sigma $%
-immanent relations (\ref{a2.3}), (\ref{a2.16b}). In the proper Euclidean
space the tube (\ref{b5.3}) degenerates into the straight line, passing
through the point $Q_0$ collinear to the given vector ${\bf P}_0{\bf P}_1$.

Let us define the set $\omega _{Q_0}=\left\{ {\bf Q}_0{\bf Q)|}Q\in \Omega
\right\} $ of vectors ${\bf Q}_0{\bf Q}$. Then
\begin{equation}
{\cal C}(P_0,P_1;Q_0)=\left\{ {\bf Q}_0{\bf Q|}Q\in {\cal T}%
(P_0,P_1;Q_0)\right\} \subset \omega _{Q_0}  \label{b5.4}
\end{equation}
is the collinearity cone of vectors ${\bf Q}_0{\bf Q}$ collinear to vector $%
{\bf P}_0{\bf P}_1$. Thus, the one-dimensionality of the first order tubes
and the collinearity cone degeneration are connected phenomena.

In the Riemannian geometry the very special property of the proper Euclidean
geometry (the collinearity cone degeneration) is considered to be a property
of any geometry and extended to the case of Riemannian geometry. The line $%
{\cal L}$, defined as a continuous mapping (\ref{b5.5}) is considered to be
the most important geometric object. This object is considered to be more
important, than the metric, and metric in the Riemannian geometry is defined
in terms of the shortest lines. Use of line as a basic concept of geometry
is inadequate for description of geometry and poses problems, which appears
to be artificial.

First, extension of curves introduces nonlocal features in the geometry
description. Nonlocalilty of description manifests itself: (1) in violation
of isometrical embeddability of nonconvex regions in the space, from which
they are cut, (2) in violation of absolute parallelism of vectors at
different points of space. These unnatural properties of Riemannian geometry
are corollaries of the metric definition via concept of a curve. In $\sigma $%
-Riemannian geometry such properties of Euclidean geometry as absolute
parallelism and isometrical embeddability of nonconvex regions conserve
completely. All this is a manifestation of negation of nondegeneracy. as a
natural property of geometry. But one fails to remove nondegeneracy  of
non-uniform geometry. It exists for spacelike vectors even in the Minkowski
geometry.

As far as one cannot remove nondegeneracy from Riemannian geometry, it
seems reasonable to recognize that the nondegeneracy is a natural geometric
property, and T-geometric conception is more perfect, than the Riemannian
conception of geometry. A corollary of this conclusion is a reconstruction
of local description and absolute parallelism (the last may be useful for
formulation of integral conservation laws in a curved space-time). Besides,
the T-geometric conception is essentially simpler, than the Riemannian one.
It has simpler structure and uses simpler method of description. The
fundamental mapping (\ref{a2.1dc}), introducing multivector in T-geometry is
essentially simpler than fundamental mapping (\ref{b5.5}), introducing the
curve in Riemannian geometry. The mapping (\ref{a2.1dc}) deals with finite
objects. It does not contain any references to limiting processes, or limits
whose existence should be provided. Finally, the T-geometry is not sensitive
to that, whether the real space-time is continuous, or only fine-grained.
This is important also, because it seems not to be tested experimentally.

\end{document}